\documentclass[letterpaper]{report}
\usepackage{graphics}
\usepackage{amsfonts}
\usepackage{amssymb}
\usepackage[mathcal]{eucal}
\usepackage{mathrsfs}
\usepackage{latexsym}
\usepackage{amsthm}
\usepackage{amsmath,amscd}
\usepackage{enumerate}
\usepackage{epic}
\usepackage[letterpaper]{hyperref}

\renewcommand{\P}{\mathscr{P}}

\newcommand{\G}{\mathscr{G}}
\renewcommand{\H}{\mathscr{H}}
\newcommand{\B}{\mathscr{B}}
\newcommand{\T}{\ensuremath{\mathscr{T}}}
\newtheorem{theorem}{Theorem}[section]
\newtheorem{lemma}[theorem]{Lemma}

\newtheorem{proposition}[theorem]{Proposition}
\newtheorem{corollary}[theorem]{Corollary}
\newtheorem{definition}[theorem]{Definition}

\newcommand{\leftexp}[2]{{\vphantom{#2}}^{#1}{#2}}
\newcommand{\sub}{\subseteq}
\newcommand{\sm}{\smallsetminus}

\newcommand{\rng}{\mathop{\mathrm{rng}}}

\newcommand{\fct}{\mathop{\longrightarrow}}

\newcommand{\Ult}{\mathop{\mathrm{Ult}}}
\newcommand{\Clop}{\mathop{\mathrm{Clop}}}
\newcommand{\cmpn}{\mathop{\mathrm{cmpn}}}
\newcommand{\res}{\upharpoonright}
\renewcommand{\phi}{\varphi}

\newcommand{\FR}[1]{\mathop{\mathrm{Fr}}\left(#1\right)}
\newcommand{\FinCo}[1]{\mathop{\mathrm{FinCo}}\left(#1\right)}
\newcommand{\At}[1]{\mathop{\mathrm{At}}\left(#1\right)}

\newcommand{\BC}[1]{\mathscr B_c\left(#1\right)}
\newcommand{\BA}[1]{\mathscr B_a\left(#1\right)}
\newcommand{\PFBA}{\ensuremath{\bot{\mathrm f}\mathcal B\mathcal A} }

\newcommand{\bas}{boolean algebras}

\newcommand{\ba}{boolean algebra}
\newcommand{\CSP}{countable separation property}
\newcommand{\Wolog}{Without loss of generality}
\newcommand{\wolog}{without loss of generality}
\newcommand{\PInd}{n\mathrm{Ind}}

\newcommand{\pa}[1]{\left(#1\right)}
\newcommand{\br}[1]{\left[#1\right]}
\newcommand{\Bco}[1]{\mathscr B_{co}\pa{#1}}
\newcommand{\Baco}[1]{\mathscr B_{aco}\pa{#1}}
\newcommand{\J}{\mathfrak i}
\newcommand{\s}[1]{\left\{#1\right\}}
\newcommand{\pc}[1]{\ensuremath{\left(\perp #1\right)}}
\newcommand{\defeq}{\ensuremath{\stackrel{\text{\tiny def}}{=}}}
\renewcommand{\epsilon}{\varepsilon}
\renewcommand{\emptyset}{\varnothing}
\newcommand{\nor}[1]{\left\|#1\right\|}
\renewcommand{\bar}{\overline}
\newcommand{\K}{\mathcal K}
\newcommand{\subs}[2]{\left[#1\right]^{#2}}
\newcommand{\F}{\mathcal F}
\newcommand{\FK}{\F\circ \K}

\title{Variations of independence in Boolean Algebras}

\author{Corey T.~Bruns}
\date{March 13, 2008}
\begin{document}
\maketitle

\copyright 2008 Corey T.~Bruns\\
This work is licensed under the Creative Commons Attribution-Noncommercial-Share Alike 3.0 United States License. To view a copy of this license, visit

  \url{http://creativecommons.org/licenses/by-nc-sa/3.0/us/} 

or send a letter to Creative Commons, 171 Second Street, Suite 300, San Francisco, California, 94105, US

\tableofcontents
\listoffigures
\chapter{$n$-free Boolean Algebras}
\section{Definition of $n$-free}
A \ba{}  $A$ is free over its subset $X$ if it has the universal property that every function $f$ from $X$ to a \ba{}  $B$ extends to a unique homomorphism.  This is equivalent to requiring that $X$ be independent and generate $A$ (uniqueness).   A generalization, $\perp$-free, is introduced in Heindorf \cite{Heindorf94}, and some of its properties are dealt with.  I follow his notation for some of its properties, but that of Koppelberg \cite{Handbook} for the operations $+,\cdot,-,0,1$ on Boolean Algebras, with the addition that for an element $a$ of a \ba{}, we let $a^0=-a$ and $a^1=a$.  An elementary product of $X$ is an element of the form $\prod_{x\in R}x^{\epsilon_x}$ where $R$ is a finite subset of $X$ and $\epsilon\in \leftexp R2$.  We further generalize the notion of freeness to $n$-freeness for $1\leq n \leq \omega$.

It is nice to have a symbol for disjointness; we define $a\perp b $ if and only if $a\cdot b =0$.

\begin{definition}\label{def:npp} Let $n$ be a positive integer, $A$ and $B$ be nontrivial \bas{}, and $U\sub A$.  A function
$f:U\rightarrow B$ is $n$-preserving if and only if for every
$a_0,a_1,\ldots,a_{n-1}\in U$, $\prod_{i<n}a_i=0$ implies that
$\prod_{i<n}f\pa{a_i}=0$.
\end{definition}

We may also do this for $\omega$.
\begin{definition}\label{def:opp} Let $A$ and $B$ be nontrivial \bas{}, and $U\sub A$.  A function
$f:U\rightarrow B$ is $\omega$-preserving if and only if for every
finite $H\sub U$, $\prod H=0$ implies that
$\prod f\left[H\right]=0$.
\end{definition}

Then we say that $A$ is $n$-free over $X$ if every $n$-preserving function from $X$ into arbitrary $B$ extends to a unique homomorphism.  The uniqueness just requires that $X$ be a generating set for $A$: 
\begin{proposition}Let $X\sub A^+$ be such that every $n$-preserving function $f:X\rightarrow B$ extends to a homomorphism. If every such $f$ extends uniquely, then $X$ generates $A$.
\end{proposition}
\begin{proof}We show the contrapositive.

Suppose that $X$ does not generate $A$, that is, there is a
$y\in A\sm\left<X\right>$. \Wolog, $A=\left<X\right>\pa y$. Let $B=\s{0,1}$ and $f$ be the zero function on
$X$.  We will show that $g\defeq f\cup\s{\pa{y,0}}$ and $h\defeq f\cup\s{\pa{y,1}}$ both extend
to homomorphisms from $A$ to $B$, thus $f$ does not extend uniquely.

We use Sikorski's extension criterion--theorem 5.5 of Koppelberg \cite{Handbook}.

Let $R\in\subs X{<\omega}, \epsilon\in \leftexp R 2$, and $\delta\in 2$ be such that $\prod_{x\in R}x^{\epsilon_x}\cdot y^\delta=0$.
If $\delta = 1$, then $\prod_{x\in R}g\pa x ^{\epsilon_x}\cdot g\pa y = \prod_{x\in R}g\pa x ^{\epsilon_x}\cdot 0=0$, as required.  If $\delta=0$, then $y\leq \prod_{x\in R}x^{\epsilon_x}$ and $g\pa y =0\leq\prod_{x\in R}g\pa x^{\epsilon_x}$, as required, so that $g$ extends to a homomorphism.

If $\delta = 0$, then $\prod_{x\in R}h\pa x ^{\epsilon_x}\cdot h\pa y = \prod_{x\in R}h\pa x ^{\epsilon_x}\cdot 0=0$, as required.  If $\delta=1$, then $y\leq \prod_{x\in R}x^{\epsilon_x}$ and $h\pa y =0\leq\prod_{x\in R}h\pa x^{\epsilon_x}$, as required, so that $h$ extends to a homomorphism.
\end{proof}

The existence of such extensions is equivalent to an algebraic property of $X$, namely that $X^+$ is $n$-independent. This notion is defined below, and the equivalence is proved.  (For $n=1$, this is the usual notion of free and independent; for $n=2$, the notions are called $\perp$-free and $\perp$-independent by Heindorf \cite{Heindorf94};  Theorem 1.3 in the same paper shows that a $2$-free \ba{} has a $2$-independent generating set.  We differ from Heindorf in that he allows $0$ to be an element of a $\perp$-independent set.) Clearly any function that is $n$-preserving is also $m$-preserving for all $m\leq n\leq \omega$, so that an $m$-free \ba{} is also $n$-free over the same set; in particular, any $n$-free \ba{} is $\omega$-free.

Freeness over $X$ with Sikorski's extension criterion implies that no elementary products over $X$ can be $0$.  $n$-independence weakens this by allowing products of $n$ or fewer elements of $X$ to be $0$.  This requires some other elementary products to be $0$ as well--if $x_1\cdot x_2\cdot\ldots\cdot x_m=0$, then any elementary product that includes $x_1,\ldots,x_m$ each with exponent $1$ must also be $0$.

\begin{definition} Let $A$ be a \ba{}.  For $n$ a positive integer, $X\sub A$ is $n$-independent if and only if $0\notin X$ and for all nonempty finite
subsets $F$ and $G$ of $X$, the following three conditions hold:
\begin{description}
\item[$\left(\perp 1\right)$] $\sum F \neq 1.$
\item[$\left(\perp 2\right)_n$] If $\prod F =0$, there is an $F'\sub F$ with
$\left|F'\right|\leq n$ such that $\prod F'=0$.
\item[$\left(\perp 3\right)$] If $0\neq \prod F \leq \sum G$, then $F\cap
G\neq\emptyset$.
\end{description}
\end{definition}

\begin{definition} Let $A$ be a \ba{}.  $X\sub A$ is $\omega$-independent if and only if $0\notin X$ and for all nonempty finite
subsets $F$ and $G$ of $X$, the following two conditions hold:
\begin{description}
\item[$\left(\perp 1\right)$] $\sum F \neq 1.$
\item[$\left(\perp 3\right)$] If $0\neq \prod F \leq \sum G$, then $F\cap
G\neq\emptyset$.
\end{description}
\end{definition}

We note that in both the above definitions, if $X$ is infinite, then
$\pc3\Rightarrow\pc1$; suppose \pc1 fails; take a finite $G$ with $\sum G=1$, then take some $x\notin G$
and let $F\defeq\s x$; then $0<\prod F\leq\sum G$ and $F\cap G\neq \emptyset$.

\pc3 has several equivalent forms which will be useful in the sequel.
\begin{proposition}\label{prop:a-2}The following are equivalent for a subset of $X$ of a \ba{}
A:
\begin{enumerate}
\item For all nonempty finite $F,G\sub X$, \pc3.
\item For all nonempty finite $F,G\sub X$ such that $F\cap G=\emptyset$ and
$\prod F\neq 0$, $\prod F \not\leq \sum G$.
\item For all nonempty finite $F,G\sub X$ such that $F\cap G=\emptyset$ and
$\prod F\neq 0$, $\prod F \cdot \prod -G\neq 0$, where $-G\defeq \s{-g: g\in
G}$.
\item Let $X$ be bijectively enumerated by $I$ such that $X=\s{x_i:i\in I}$. For all nonempty finite $R\sub I$ and all $\epsilon\in \leftexp R 2$ such
that $1\in\rng\epsilon$ and $\prod_{\stackrel{i\in R}{\epsilon_i = 1}}x_i\neq 0$,
$\prod_{i\in R}x_i^{\epsilon_i}\neq 0$. 
\end{enumerate}
\end{proposition}
In words, the final equivalent says that no elementary product of elements of $X$ is
$0$ unless the product of the non-complemented elements is $0$.  We note that
in the presence of $\pc2_n$, the words ``of $n$'' may be inserted after
``product.''
\begin{proof}
We begin by pointing out that \pc3 has two hypotheses, $0\neq \prod F$ and
$\prod F\leq\sum G$.  Thus the contrapositive of \pc3 is ``If $F\cap
G=\emptyset$, then $0=\prod F$ or $\prod F\not\leq\sum G$,'' which is
equivalent to (2).

(2) and (3) are equivalent by some elementary facts: $a\leq b \;\iff a\cdot-b=0$ and de Morgan's law that $-\sum G=\prod -G$.

(3) $\Rightarrow$ (4):

Assume (3) and the hypotheses of (4).  If $\rng \epsilon =\s1$, the conclusion is clear.  Otherwise, let $F\defeq\s{x_i:i\in R\mbox{ and } \epsilon_i=0}$.  Then (3) implies that $\prod_{i\in R}x_i^{\epsilon_i}\neq 0$, as we wanted.

(4) $\Rightarrow$ (3):

Assume (4) and the hypotheses of (3).  Let $R\defeq\s{i\in I : x_i\in F\cup G}$ and let $\epsilon_i=1$ if $x_i\in F$ and $\epsilon_i=0$ otherwise.  Then (4) implies that $\prod F\cdot\prod -G \neq 0$, as we wanted.
\end{proof}

\begin{proposition}\label{prop:a-1}The following are equivalent for a subset $X$ of a \ba{}
A:
\begin{enumerate}

\item $X$ is $\omega$-independent
\item Let $X$ be bijectively enumerated by $I$ such that $X=\s{x_i:i\in I}$. For all nonempty finite $R\sub I$ and all $\epsilon\in \leftexp R 2$ such
that $\prod_{\stackrel{i\in R}{\epsilon_i = 1}}x_i\neq 0$,
$\prod_{i\in R}x_i^{\epsilon_i}\neq 0$.
\end{enumerate}
\end{proposition}
\begin{proof}

The proof is the similar to that of proposition \ref{prop:a-2}.  \pc1 is taken care of since products over an empty index set are taken to be $1$ by definition.
\end{proof}

In the same spirit, we have an equivalent definition of $n$-independent.
\begin{proposition}\label{prop:a}Let $n$ be a positive integer or $\omega$, $A$ a nontrivial \ba{} and $X\sub A^+$.  $X$ is $n$-independent if and only if for every $R\in\left[X\right]^{<\omega}$ and every $\epsilon\in\leftexp R2$, if $\prod_{x\in R}x^{\epsilon_x}=0$ then there is an $R'\sub R$ with $\left|R'\right|\leq n$ such that $\epsilon\left[R'\right]=\s1$ and $\prod R' =0$.
\end{proposition}
\begin{proof}
If $n=\omega$, this is part of proposition \ref{prop:a-1}.

Let $n$ be a positive integer, $A$ a \ba{}, and $X\sub A^+$.

We first show that $n$-independent sets have the indicated property.

Assume that $X$ is $n$-independent; take $R\in\left[X\right]^{<\omega}$ and $\epsilon\in\leftexp R2$ such that  $\prod_{x\in R}x^{\epsilon_x}=0$.  Let $F=\s{x\in R:\epsilon_x=1}$ and $G=\s{x\in R:\epsilon_x=0}$.  $F\neq\emptyset$; otherwise $\sum R=-\prod R = 1$, contradicting \pc1.  Since $\prod_{x\in R}x^{\epsilon_x}=\prod F \cdot \prod -G$, we have that $\prod F\leq\sum G$. 
If $G=\emptyset$, then $\sum G=0$ and so $\prod F =0$ as well.  If $G\neq \emptyset$, then $\prod F=0$ since $F\cap G=\emptyset$, using \pc3.  Then $R'$ is found by $\pc2_n$.
Now we show that sets with the indicated property are $n$-independent.

Assume that $X$ has the indicated condition and $F,G\in\left[X\right]^{<\omega}\sm\s\emptyset$.  We have three conditions to check.
\begin{description}
\item[\pc1]  Suppose that $\sum F=1$.  We let $F$ be the set $R$ in the condition, setting $\epsilon_x=0$ for all $x\in F$.  Then $\prod_{x\in F}x^{\epsilon_x}=\prod -F=-\sum F = 0$ and $\s{x\in F:\epsilon_x=1}=\emptyset$, thus there is no $R'$ as in the condition, since products over an empty index set are equal to 1.

\item[$\pc2_n$] Suppose that $\prod F=0$.  Again we let $F$ be the set $R$ in the condition, this time setting $\epsilon_x=1$ for all $x\in F$.  Then the condition gives us the necessary $F'$

\item[\pc3]  Suppose that $0\neq\prod F\leq\sum G$ and $F\cap G=\emptyset$.  Let $R=F\cup G$ and $\epsilon\in \leftexp R2$ be such that $\epsilon\left[F\right]=\s1$ and $\epsilon\left[G\right]=\s0$.  Then $\prod_{x\in R}x^{\epsilon_x}=0$ and the condition gives $\prod F=0$, which contradicts the original supposition.
\end{description}
\end{proof}

$n$-independence is a strictly weaker property than independence.  It is not hard to construct sets which are $n$-independent but not $\pa{n-1}$-independent, however it is more difficult to show that there is an algebra which is $n$-free and not $\pa{n-1}$-free.  This will be done for infinitely many $n$ later.

\begin{lemma}\label{pindtoind} If $H$ is an 
$\omega$-independent set that has no finite subset $F$ such that $\prod F=0$, $H$ is in fact
independent.  Furthermore, if $H$ is $n$-independent with no subset $F$ of size $n$ or less with $\prod F=0$, then $H$ is independent.\end{lemma}
\begin{proof}
We only need show that $\pc2_1$ holds, which it does vacuously.
\end{proof}

$2$-independence, and thus $n$-independence for $2\leq n\leq\omega$, is also a generalization of pairwise disjointness on
infinite sets.

\begin{theorem}\label{ipd}If $X\sub \B^+$ is an infinite pairwise disjoint set, then $X$ is
$2$-independent.\end{theorem}
\begin{proof}

This is clear from proposition \ref{prop:a}.
\end{proof}

Some non-trivial examples of $2$-free \bas{} are the finite-cofinite algebras.  For infinite $\kappa$, let $A=\mathop{\mathrm FinCo}\left(\kappa\right)$.  $\mathrm{At}\pa A$  is a $2$-independent generating set for $A$.  That $\mathrm{At}\pa A$  is a generating set is clear, and it is $2$-independent by theorem \ref{ipd}.

The $n=2$ cases of the preceding lemma and theorems, though niether explicit nor implicit in Heindorf \cite{Heindorf94}, were likely (in our opinion) to have been motivation for the definition.

Having an
$n$-independent generating set is equivalent to $n$-freeness. This is known in
Koppelberg \cite{Handbook} for $n=1$ and Heindorf \cite{Heindorf94} for $n=2$.  Our proof is more elementary than that of Heindorf \cite{Heindorf94} in that it avoids clone theory.

\begin{theorem}\label{thm:ofree_oindep}If $A$ is $\omega$-free over $X$, then
$X^+$ is $\omega$-independent.
\end{theorem}
\begin{proof}
Let $A$ and $X$ be as in the hypothesis; we show that $X^+$ is
$\omega$-independent.

\Wolog{}, we may assume that $0\notin X$ so that $X^+=X$.
\begin{description}
\item[\pc1]Let $f:X\rightarrow \s{0,1}$ be such that $f\br X =\s 0$.  Clearly $f$ is
$\omega$-preserving and thus extends to a homomorphism $\overline f$.  Take
$F\in\left[X\right]^{<\omega}$; then $\overline f \pa{\sum F}=\sum
f\left[F\right]=0$, so that $\sum F\neq 1$.
\item[\pc3]Take $F,G\in\left[X\right]^{<\omega}$ such that $F\cap G=0$ and
$\prod F\neq 0$.  Let $f:X\rightarrow\s{0,1}$ be such that
$f\left[F\right]=\s1$ and $f\left[X\sm F\right]=\s0$.  We claim that $f$ is
$\omega$-preserving.  If $H\sub X$ is finite such that
$\prod f\left[H\right]\neq0$, then it must be
that $H\sub F$, and hence $\prod H\neq 0$.  Thus $f$ extends to a homomorphism $\overline
f$.  Then $\overline f\pa{\prod F\cdot \prod -G}=\prod
f\left[F\right]\cdot\prod\overline f\left[-G\right]=1$, and so $\prod
F\cdot\prod -G\neq 0$.
\end{description}

\end{proof}

\begin{theorem}\label{thm:nfree_nindep}  Let $n$ be a positive integer and $A$ a \ba{}.  If $A$ is $n$-free over $X$, then $X^+$ is $n$-independent.
\end{theorem}

\begin{proof} Again, \wolog{} $X=X^+$.

From theorem \ref{thm:ofree_oindep}, $X$ is $\omega$-independent, so we need
only show that $\pc2_n$ holds for $X$.  We do this by contradiction; assume
that $F\sub X$ is finite, of cardinality greater than $n$, $\prod F=0$, and
every subset $F'\sub F$ where $F'$ is of size $n$ is such that $\prod F'\neq
0$.

Define $f:X\rightarrow \s{0,1}$ by letting $f\left[F\right]=\s1$ and
$f\left[X\sm F\right]=\s0$.

Then $f$ is $n$-preserving. Let $G\sub X$ be of size $n$ and have $\prod G=0$.
Then $G\not\sub F$, so some $x\in G$ has $f\pa x =0$, so  $\prod f\left[G\right] = 0$.  Thus $f$ must extend to a homomorphism, but then $f\pa 0 =f\pa{\prod F}=\prod f\left[F\right]=\prod \s1=1$, which is a contradiction.
\end{proof}

\begin{theorem}\label{thm:oindep_ofree}Let $A$ be generated by its
$\omega$-independent subset $X$.  Then $A$ is $\omega$-free over $X$.
\end{theorem}
\begin{proof}
This is an easy application of Sikorski's extension criterion.  Let $f$ be
an $\omega$-preserving function with domain $X$; we will show that
$f$ extends to a unique homomorphism.

Take a finite $H\sub X$ and $\epsilon\in\leftexp H2$ such that $\prod_{h\in
H}h^{\epsilon_h}=0$.  Then by \pc3 and \pc1, $\prod_{\epsilon_h=1}h=0$.  Then
since $f$ is $\omega$-preserving, $\prod_{\epsilon_h=1}f\pa h =0$ and
thus $\prod_{h\in H}f\pa h ^{\epsilon_h}=0$.  Thus by Sikorski's extension criterion, $f$ extends to a homomorphism. 

Uniqueness is clear as $X$ is a generating set.
\end{proof}

\begin{theorem}\label{thm:nindep_nfree}Let $n$ be a positive integer.  If $X$ generates $A$ and $X$ is $n$-independent, then $A$ is
$n$-free over $X$
\end{theorem}
\begin{proof}
This is also a straightforward application of Sikorski's extension criterion.  Let
$f$ be an $n$-preserving function with domain $X$.  

Take any distinct $x_0,x_1,\ldots,x_{k-1}$ and $\epsilon\in\leftexp n2$ such that
$\prod_{i<k}x_i^{\epsilon_i}=0$.  

Then by proposition \ref{prop:a}, there is an $F'\sub\s{x_i:\epsilon_i=1\;\mbox{and}\; i<k}$ such that $\left|F'\right|\leq
n$ and $\prod F'=0$.  Since $f$ is $n$-preserving, it must be that
$\prod f\left[F'\right]=0$, and thus $\prod_{i<k}f\pa{x_i}^{\epsilon_i}=0$.
Thus, by Sikorski's extension criterion, $f$ extends to a homomorphism. 

Uniqueness is clear as $X$ is a generating set.
\end{proof}

So we have shown that the universal algebraic property defining $n$-free \bas{}
is equivalent to having an $n$-independent generating set.

\begin{theorem}\label{semigroup} $\omega$-free \bas{} (and thus all $n$-free
\bas{}) are semigroup algebras.
\end{theorem}
A semigroup algebra is a \ba{}  that has a generating set that includes $\left\{0,1\right\}$, is closed under the product operation, and is disjunctive when $0$ is removed.
\begin{proof}
Let $A$ be $\omega$-free over $G$.  Then let $H'$ be the closure of $G\cup \left\{0,1\right\}$ under finite products, that is, the set of all finite products of elements of $G$, along with $0$ and $1$.  Clearly $H'$ generates $A$, includes $\left\{0,1\right\}$, and is closed under products, so all that remains is to show that $H=H'\sm\left\{0\right\}$ is disjunctive.
From proposition 2.1 of Monk \cite{CINV2}, $H$ is disjunctive if and only if for every $M\sub H$ there is a homomorphism $f$ from $\left<H\right>$ into $\P\left(M\right)$ such that $f\left(h\right)=M\downarrow h$ for all $h\in H$. 

To this end, given $M\sub H$, let $f:G\rightarrow \P\left(M\right)$ be defined by $g\mapsto M\downarrow g$.  We claim that $f$ is $\omega$-preserving.  Suppose $G'\in \left[G\right]^{<\omega}$ is such that $\prod G'=0$.  Then $\prod_{g\in G'}f\pa g = \bigcap_{g\in G'}\left(M\downarrow g \right)=\left\{a\in M :\forall g\in G'\left[a\leq g\right] \right\}=\emptyset$.  So $f$ extends to a unique homomorphism $\hat f$ from $A$ to $\P\left(M\right)$.  If $h\in H\sm \left\{1\right\}$, then $h=g_1 \cdot g_2 \cdot \ldots \cdot g_n$ where each $g_i\in G$.  So $$\hat f\left(h\right)=\hat f\left(g_1 \cdot g_2 \cdot \ldots \cdot g_n\right)=f\left(g_1\right)\cap f\left(g_2\right)\cap \ldots \cap f\left(g_n\right)=$$ $$=\left(M\downarrow g_1 \right)\cap \left(M \downarrow g_2\right)\cap \ldots \cap \left(M\downarrow g_n\right) =M\downarrow \left(g_1 \cdot g_2 \cdot \ldots \cdot g_n\right)=M\downarrow h.$$  Likewise, $\hat f\left(1\right)=M = M\downarrow 1$.  Thus $H$ is disjunctive and $A$ is a semigroup algebra over $H'$.
\end{proof}

\section{Definition of Graph Space}\label{sec:gs}
Given a graph $\G=\left<G,E\right>$, Bell \cite{Bell82} defines a topology on
the set of complete subgraphs (cliques) of $\G$.    We will also use a similar construction on the set of anticliques.  By anticlique, we mean a set of vertices whose induced subgraph has no edges--we use this term as it is more suggestive of the relationship to cliques (a clique in $\G$ is an anticlique in $ \bar\G$ and vice versa, where $\bar\G$ is the complement graph of $\G=\left<G,E\right>$.  That is,   $\bar\G\defeq\left<G,\left[G\right]^2\sm E\right>$.) than the term ``stable set'' and less confusing in this context than ``independent set.'' All of our graphs will be loopless.  It is useful to identify subsets of $G$ with their
characteristic functions in $\leftexp G2$, so we give the definition using
this identification.

\begin{definition}For a graph $\G=\left<G,E\right>$,  let
$$C\left(\G\right)=\left\{f\in \leftexp G 2 :
f^{-1}\left[\left\{1\right\}\right] \mbox{is a complete subgraph of } \G\right\}.$$  For $v\in G$, we define two subsets of $C\left(\G\right)$, $$v_+=\left\{f\in C\left(\G\right) : f\left(v\right)=1\right\}\mbox{ and }v_-=\left\{f\in C\left(\G\right) : f\left(v\right)=0\right\}.$$  These are the subsets of $C\left(\G\right)$ corresponding to the sets of complete subgraphs of $\G$ including and excluding $v$, respectively.  Then $\bigcup_{v\in G}\left\{v_+, v_-\right\}$ is a subbase for a topology on $C\left(\G\right)$, the clique graph space topology for $\G$.  

Similarly, let $A\left(\G\right)=\left\{f\in \leftexp G 2 :
f^{-1}\left[\left\{1\right\}\right] \mbox{is an independent subgraph of } \G\right\}$.  $v_+$ and $v_-$ are defined similarly, and the topology on $A\left(\G\right)$ with subbase $\bigcup_{v\in G}\left\{v_+, v_-\right\}$ is called the anticlique graph space topology for $\G$.  
\end{definition}

Note that the clique graph space of $\G$ is the anticlique graph space of $ \bar\G$, and vice versa, so if $\G$ is left unspecified, we may omit the word (anti)clique.

The following observation (made in Bell \cite{Bell82}) is the reason for using the
characteristic functions; the proof is trivial.

\begin{proposition}
Let $2$ have the discrete topology and $\leftexp G2$ the product topology.  $C\left(\G\right)$ with the graph space topology is the subspace of $\leftexp G2$ consisting of all $f:G\rightarrow2$ such that $f\left(g\right)=f\left(h\right)=1$ implies $\left\{g,h\right\}\in E$.
Any such space (or one homeomorphic to such a space) is a graph space.

Such a space is closed in $\leftexp G2$, so is a Boolean space.

\end{proposition}
Note as in section 3 of Bell \cite{Bell82}, that for any set  $I$ the Cantor space $\leftexp I2$ is a graph space; in particular, it is $C\left(\left<I,\left[I\right]^2\right>\right)$.

The graph associated with a given clique graph space is not unique, even up to
isomorphism.  In the case of a finite graph, the clique graph space is the discrete
topology on the set of cliques, so any two graphs with the same
number of cliques have homeomorphic clique graph spaces.  For example, the
graph on three vertices with a single edge has 5 cliques: the empty
subgraph, 3 singletons, and 1 doubleton.  The graph on 4 vertices with no edges
also has 5 cliques: 1 empty and 4 singletons.

There are also non-isomorphic infinite graphs with homeomorphic graph spaces; we will give this example in section \ref{sec:ac} after the relevant theorems about free products and their relationship with graphs are discussed in that section.

\section{Previous results on Graph Spaces and $2$-free \bas}\label{sec:prevres}
The main result of Heindorf \cite{Heindorf94} is that any second countable Boolean
space is a graph space.  The dual statement (which is what is actually proved)
is that every countable \ba{}  is $2$-free.

This gives some less trivial examples of $2$-independent sets which are not independent:  Take $A$ to be any countable \ba{}  that is not isomorphic to $\FR\omega$.  Then since $A$ is $2$-free, it has a $2$-independent generating set (of size $\omega$), and since $A$ is not free, that set can not be independent.

In proposition 3.1 of Bell \cite{Bell82}, it is shown that every graph space is supercompact, that is,
that every graph space has a subbase with the property that any open cover
consisting of subbasic sets has a 2-element subcover.  The closed set dual of this is sometimes more handy--a binary set is a set $\B$ with the property that if any set $\mathcal S\sub \B$ has empty intersection, there are two elements of $\mathcal S$ that are disjoint.  Then the closed set dual is that the space has a binary subbase for the closed sets. 

We will later (theorem \ref{thm:fcf_not2free}, ff.) give an example of a supercompact boolean space that is not a graph space.

This may be used to show that not all (uncountable) \bas{} are $2$-free--$\beta\pa\omega \sm \omega$ is not supercompact (Corollary 1.1.6, of van Mill\cite{MR0464160}), so its dual algebra, which is isomorphic to $\P\pa\omega /\mbox{fin}$ is not $2$-free.  In fact, we later show (after lemma \ref{finco_is_image}) that $\P\pa\omega /\mbox{fin}$ is not $\omega$-free.

The phrase ``graph space'' has been used for a different notion in the literature. van Mill \cite{MR0464160} follows de Groot \cite{MR0372808} and defines a subbase on the set of maximal anticliques of a graph, letting the subbase be $\s{v_+ : v\in G}$.  van Mill \cite{MR0464160} calls these spaces ``graph spaces'' (de Groot calls them ``stability spaces'') and repeats de Groot's \cite{MR0372808} implicit result that these spaces are the supercompact spaces, that is, the stability space of every graph is supercompact, and every supercompact space is the stability graph of the $\perp$-graph (defined in the next section) of its binary closed subbase.
We show in theorem \ref{thm:fcf_not2free} that there is a graph space which is
not supercompact.

Kunen \cite{KunenSet} constructs some \bas{} from partial orders (which are
somewhat graph-like) in Chapter
II, section 3 and exercise C5 of chapter VII; we do not know if these \bas{} or
some reasonable subalgebras are $\omega$-free.  Galvin \cite{MR585558} also
constructs some \bas{} with unusual properties from partial orders; we do not
know if they are $\omega$-free either.

\section{Graphs and their (Anti)clique Algebras}

From a graph $\G=\left<G,E\right>$, one constructs an (anti)clique graph space and then from the graph space, a
boolean algebra.  One may skip the graph space step,  directly constructing 
boolean algebras on $\P\left(C\left(\G\right)\right)$ and $\P\left(A\left(\G\right)\right)$ generated by $\left\{v_+ : v\in
G\right\}$.   We call these \bas{} the clique algebra of $\G$ and the anticlique algebra of $\G$, denoted $\BC\G$ and $\BA\G$ respectively.   We now show that these \bas{} are $2$-free; to do so, we need only show that
$\left\{v_+ : v\in G\right\}$ is $2$-independent in $\P\left(C\left(\G\right)\right)$ and $\P\left(A\left(\G\right)\right)$  Actually, it is enough to do this for just $\P\left(C\left(\G\right)\right)$ as $\BA\G=\BC{ \bar\G}$.  It is often necessary to specify a $2$-independent generating set; when a \ba{}  is written as $\BA\G$ or $\BC\G$, this will be assumed to be $\left\{v_+ : v\in G\right\}$,  which we will  write  as $G_+$.  This is implicit in Heindorf \cite{Heindorf94}.

\begin{theorem}$G_+=\left\{v_+ : v\in G\right\}$ is $2$-independent in $\P\left(C\left(\G\right)\right)$.
\label{v+pind}
\end{theorem}
\begin{proof}

We will apply proposition \ref{prop:a}.  Suppose that $R\in\left[G\right]^{<\omega}, \epsilon\in \leftexp R2,$ and $\bigcap_{v\in R}v_+^{\epsilon_v}=\emptyset$.  By way of contradiction, suppose that $$\forall R'\in \left[R\right]^{\leq2}\pa{\epsilon\left[R'\right]\sub\s 1 \rightarrow \prod_{v\in R'}v_+\neq 0}.$$  Then $C\defeq\s{v\in R:\epsilon_v=1}$ is a clique and $C\in\bigcap_{v\in R}v_+^{\epsilon_v}$, contradiction.

\end{proof}

We will show an equivalence of categories between graphs, denoted by $\mathfrak G$ and $2$-free \bas{} denoted by \PFBA.  First we have some definitions:

\begin{definition}\label{graphdef} The objects in $\mathfrak G$ are graphs as
used above.  The morphisms are graph homomorphisms.  A graph homomorphism from
$\G_0=\left<G_0,E_0\right>$ to $\G_1=\left<G_1,E_1\right>$ is a function
$f:G_0\rightarrow G_1$ such that if $\left\{u,v\right\}\in E_0$, then
$\left\{f\left(u\right),f\left(v\right)\right\}\in E_1$, i.e. $f$ preserves
adjacency.
\end{definition}

Graph homomorphisms are defined and discussed in detail in Hell and Ne{\v{s}}et{\v{r}}il \cite{GraphHom}; colorings and many other graph-theoretic concepts can be defined as the existence or non-existence of certain graph homomorphisms.

\begin{definition}\label{PFBA}  The objects in \PFBA are pairs
$\left(B,G\right)$ where $B$ is a \ba{}  and $G\sub B^+$ such that $B$ is
$2$-free over $G$.  A morphism between $\left(B,G\right)$ and
$\left(C,H\right)$ is a $2$-preserving function between their generating
sets, i.e. an $f:G\rightarrow H$ such that for all $a,b\in G$ with $a\perp b$,
$f\left(a\right)\perp f\left(b\right)$.  Since such a map extends to a unique
BA homomorphism from $B$ to $C$, we will use the same symbol for both a
morphism and its extension when there is no possibility of confusion.
\end{definition}

For the remainder of this section, $\BA\G$ and $\BC\G$ will be  objects in \PFBA, namely the pair $\left(B,G_+\right)$ where $B$ is
the (anti)clique algebra of $\G$ and $G=\left\{v_+:v\in G\right\}$ is the usual
$2$-independent generating set.  This should be regarded as the formal
definition of $\BA\G$ and $\BC\G$.

Now with each object $\pa{B,G}$ in \PFBA, we associate a graph $\K\pa{B,G}$ as follows.  The set of vertices of  $\K\pa{B,G}$ is $G$ itself.  Let $E\sub\subs G2$ be such that $\s{h,k}\in E$ if and only if $h\perp k$.  We call $\K\pa{B,G}$ the $\perp$-graph of $\pa{B,G}$.

\begin{proposition}\label{prop:A} If $\pa{B,G}\in \PFBA$, then $\pa{B,G}$ is isomorphic to $\BA{\K\pa{B,G}}$
\end{proposition}
\begin{proof}
For any $g\in G$ let $f\pa g= g_+$.  Clearly $f$ is a bijection from $G$ to $G_+$.  $f$ is $2$-preserving: if $h\perp k$, then $\s{h,k}\in E$, so $\s{h,k}$ is not an anticlique, so $h_+\perp k_+$.  Similarly, $f^{-1}$ is $2$-preserving, so $f$ extends to the desired isomorphism.
\end{proof}

\begin{proposition}\label{prop:B}If $f:\pa{B,G}\rightarrow\pa{B',G'}$ is a morphism, then $\K\pa f$ is a morphism from $\K\pa{B,G}$ to $\K\pa{B',G'}$.\end{proposition}
\begin{proof}If $\s{h,k}$ is an edge of $\K\pa{B,G}$, then $h\perp k$; thus $\s{f\pa h , f\pa k}$ is an edge of $\K\pa{B',G'}$.
\end{proof}

Now we begin to define a functor:  for each graph $\G$, let $\F\pa\G \defeq \BA\G$.

\begin{proposition}\label{prop:C}For any graph $\G$, the $\perp$-graph $\G'$ of $\F\pa\G$ is isomorphic to $\G$.
\end{proposition}
\begin{proof}
For each $g\in G$ let $f\pa g = g_+$.  Clearly $f$ is a bijection from $G$ to $G_+$.  Moreover, $\s{h,k}$ is an edge of $\G$ if and only if $h_+\perp k_+$ if and only if $\s{h_+,k_+}$ is an edge of $\G'$.
\end{proof}

We continue to define the functor $\F$; if $f:\G\rightarrow\G'$ is a morphism, define $\F\pa f : \BA\G\rightarrow \BA{\G'}$ by setting, for any $g\in G$, $\F\pa f \pa{g_+}=\pa{f\pa g}_+$.

\begin{proposition}\label{prop:D}If $f:\G\rightarrow \G'$ is a morphism, then $\F\pa f :\BA\G\rightarrow\BA{\G'}$ is a morphism.
\end{proposition}
\begin{proof} If $h_+\perp k_+$, then $\s{h,k}$ is an edge of $\G$, hence $\s{f\pa h , f\pa k}$ is an edge of $\G'$, so $\pa{f\pa h}_+\perp \pa{f\pa k}_+$.
\end{proof}

\begin{proposition}\label{prop:E}The categories \PFBA and $\mathfrak G$ are equivalent.
\end{proposition}
\begin{proof}
We follow the definitions in MacLane \cite{MR1712872}, pages 16-18.  Let $ID_\mathfrak G$ and $ID_{\PFBA}$ be the identity functors on $\mathfrak G$ and \PFBA respectively.  For a natural isomorphism of $ID_\mathfrak G$ to $\K\circ\F$, we let $f_\mathcal G$ be as defined in the proof of proposition \ref{prop:C}.  The diagram for natural transformations is then as follows:
$$\begin{CD} \G @>f_\G>>\pa{\K\circ \F}\pa\G \\
@VsVV @VV\pa{\K\circ \F}\pa sV\\
\G' @>>f_{\G'}> \pa{\K\circ \F}\pa{\G'}
\end{CD}$$

The diagram commutes: $\pa{\K\pa{\F\pa s}}\pa{f_\G\pa g}=\pa{\K\pa{\F\pa s}}\pa{g_+} = \pa{\F\pa s}\pa{g_+} = \pa{s\pa g}_+$ and $f_{\G'}\pa{s\pa g}=\pa{s\pa g}_+$.

For a natural transformation of $ID_{\PFBA}$ to $\F\circ \K$, we let $k_{\pa{B,G}}$ be defined as in the proof of proposition \ref{prop:A}. The diagram for natural transformations is then
$$\begin{CD}\pa{B,G} @>k_{\pa{B,G}}>>  \pa{\FK}\pa{B,G}\\
@VsVV  @VV\pa{\FK}\pa sV\\
\pa{B',G'} @>>k_{\pa{B',G'}}>  \pa{\FK}\pa{B',G'} 
\end{CD}$$

The diagram commutes: $\pa{\pa{\FK}\pa s}\pa{k_{\pa{B,G}} \pa g}=$\\$\pa{\F\pa{\K \pa s}}\pa{g_+}=\pa{\pa{\K\pa s }\pa g }_+=\pa{s\pa g}_+$ and $k_{\pa{B',G'}}\pa{s\pa g}=\pa{s\pa g}_+$.

\end{proof}

\begin{proposition}\label{ind}TFAE:
\begin{enumerate}
\item $H$ is a clique in $\G$.
\item $H$ is an anticlique in $ \bar\G$
\item $H_+$ is independent in $\BC\G$ and $\BA{ \bar\G}$.
\item $H_+$ is pairwise disjoint in $\BA\G$ and $\BC{ \bar\G}$.
\end{enumerate}
       
\end{proposition}
\begin{proof}

$(1) \iff (2)$ is obvious.

$(1) \Rightarrow (3)$:
If $H$ is a clique in $\G$, then $H_+$ has no disjoint pair in $\BC\G$, and is a
subset of $G_+$, so is $2$-independent.  So by \ref{pindtoind}, $H_+$ is
independent.

$(3) \Rightarrow (1)$:
If $H_+$ is independent in $\BC\G$, then every pair of vertices from $H$ is a
clique, so $H$ is a clique.

$(2)\iff (4)$:
$H$ is an anticlique in $ \bar\G$ if and only if for every $h,k\in H$, $\left\{h,k\right\}$ is not
an edge if and only if $h_+\perp k_+$ in $\BC{ \bar\G}$.

\end{proof}

These results indicate that an anticlique (independent set of vertices) in $\G$ corresponds to an independent set in $\BA\G$ and a clique corresponds to a pairwise disjoint set in $\BA\G$, We also have the corresponding cardinal function results:  if $\G$ has an anticlique of size $\kappa$, then $\mbox{Ind}\left(\BA\G\right)\geq\kappa$ and if $\G$ has a clique of size $\lambda$, then $c\left(\BA\G\right)\geq\lambda$.

We give some examples of $2$-free \bas{} with unusual properties.

For a $2$-free algebra of the form $\BC T$ for a tree (in the graph-theoretical
sense--a connected acyclic graph) or a forest $T$ of size
$\kappa$, there
are further conclusions that can be drawn.  As a forest has no triangles, all its
cliques are of size at most 2.

So any subset
of $T_+$  of size 3 or more has a disjoint pair.  

If $T$ is a $\kappa$-tree (in the order theoretic sense, 
that is, of height $\kappa$ and each level of size $<\kappa$),and we take the edge set to consist of pairs $\s{u,v}$ where $v$ is an immediate successor of $u$, then $T_+$ has
a pairwise disjoint subset of size $\kappa$--take an element of every other
level--so that $\FinCo\kappa\leq\BC T$, and $\FR\kappa\leq \BA T$. 

 It seems
difficult to avoid one of $\FinCo\kappa$ and $\FR\kappa$ as a subalgebra, as it
is necessary to find a graph of size $\kappa$ with no clique or anticlique of
size $\kappa$.  A witness to $\kappa\not\longrightarrow \pa\kappa^2_2$ is the edge set of such a graph, but we do not know about the variety of such witnesses.  If $\kappa$ is weakly compact, then there are no such witnesses and so for any graph of size $\kappa$, $\FinCo\kappa$ or $\FR\kappa$ is a subalgebra of $\BC\G$ 

As a graph can be characterized as a symmetric non-reflexive relation, for any non-reflexive relation $R$,
we may form the algebras $\BA{R\cup R^{-1}}$ and $\BC{R\cup R^{-1}}$.  When $R$
is an ordering of some sort, $R\cup R^{-1}$ is usually called the (edge set of
the) comparability
graph of $R$.  Thus for a (non-reflexive
) ordering $\left<P,<\right>$, it has comparability
graph $\G_P=\left<P,<\cup<^{-1}\right>$ and  we define its
comparability algebra $\Bco P \defeq\BC{\G_P}$ and its
incomparability algebra $\Baco P \defeq\BA{\G_P}$.  Since points in the partial order are vertices of the comparability graph, we may use the $p_+$ notation without fear of confusion. 
When $P$ is a partial order in the strict sense, $C\sub P$ is a clique in
$\G_P$ if and only if $C$ is a chain in $\leq$ if and only if $C_+$ is an independent subset of
$\Bco P$, and $A\sub P$ is an anticlique in
$\G_P$ if and only if $A$ is an antichain in $\leq$ if and only if $A_+$ is a pairwise disjoint
set in $\Bco P$.
So if $\left<T,\leq\right>$ is a $\kappa$-Suslin tree, in both $\Bco T$ and $\Baco T$, $T_+$ is a
$2$-independent set of size
$\kappa$, but has no independent subset of size $\kappa$, nor a pairwise
disjoint subset of size $\kappa$ since $T$ has neither chains nor antichains of
size $\kappa$.

\begin{proposition}
If $f:P\rightarrow Q$ is a strictly order-preserving function, that is, a morphism in the category of partial orders, then there is a homomorphism $f^*:\Baco P \rightarrow \Baco Q$ such that $f^*\pa{p_+}=f\pa p$.
\end{proposition}
\begin{proof}
By the universal property of $2$-free \bas{}, we need only show that $g$ is $2$-preserving where $g\pa{p_+}=f\pa p$; then $g$ extends to the $f^*$ of the conclusion.

Fix distinct $p,p'\in P$; if $p_+ \perp p'_+$ in $\Baco P$, then $p$ and $p'$ are comparable in $P$, without loss of generality, $p<p'$.  Then $f\pa p < f\pa{ p'}$, so that $f\pa p \perp f \pa {p'}$.
\end{proof}
Similarly, an incomparability-preserving map from $P$ to $Q$ gives rise to a homomorphism of $\Bco P$ and $\Bco Q$.

\section{Hypergraphs and their Anticlique Algebras}
There is a correspondence with hypergraphs for $\omega$-free \bas{}.
We recall that a hypergraph is a pair $\G=\left<V,E\right>$ where $V$ is called the vertex set, and $E\sub\P\pa V\sm\s\emptyset$ is called the hyperedge set; an element of $E$ is called a hyperedge. Again, we will insist on loopless hypergraphs, that is, $E\sub\P\pa V\sm\subs{V}{\leq 1}$.
 A hypergraph is $n$-uniform if $E\sub\left[V\right]^n$.  For a given hypergraph, we call a set $A\sub V$ an anticlique if it includes no hyperedges; that is, for all $e\in E$, $e\sm A\neq \emptyset$, and call the set of anticliques $A\pa \G$.

Given a hypergraph $\G$, we define an $\omega$-free \ba{} as a subalgebra of $\P\pa{A\pa \G}$.  For $v\in V$, let $v_+\defeq\s{A\in A\pa\G : v\in A}$, which is an element of $P\pa{A\pa \G}$, and for a set $H$ of vertices, $H_+\defeq\s{v_+:v\in H}$.  We then define the anticlique algebra of $\G$ as $\BA\G\defeq\left<V_+\right>$.

We do not consider cliques in general hypergraphs; there is not a unique way to define them.  For an $n$-uniform hypergraph, a clique may be succinctly defined as a set $C$ where $\subs C n \sub E$, but for a  hypergraph with hyperedges of different cardinalities, it is not clear how many hyperedges must be included in a clique.    This difficulty stems from a lack of a reasonable way to define ``complement hypergraph.''  A few possibilities for the hyperedge set of $\bar\G$ are $\P\pa G \sm E$,  $\subs G {<\omega} \sm E$,  and $\subs G {\leq\pa{\sup_{e\in E}\pa{\left|e\right|}}}\sm E$.  For an $n$-uniform hypergraph $\left<G,E\right>$, the complementary hypergraph is $\left<G,\subs G n \sm E\right>$, and then a clique in $\G$ is an anticlique in $\bar\G$.  Each possible definition for complement hypergraph results in a different definition for clique, all of which are more complicated than our definition of anticlique. Since anticliques suffice for our study, we do not choose a side on what a clique ought to be.

\begin{theorem}For any hypergraph $\G=\left<V,E\right>$, $\BA\G$ is $\omega$-free over $V_+$.\end{theorem}
\begin{proof} We need only show that $V_+$ is $\omega$-independent, we will use proposition \ref{prop:a}.

Suppose that $R\in\subs V{<\omega}, \epsilon\in \leftexp R 2$, and $\bigcap_{v\in R} v_+^{\epsilon_v}=0$.  Let $S=\s{v\in R:\epsilon_v=1}$.  If $\bigcap_{v\in S} v_+\neq 0$, let $T$ be a member of $\bigcap_{v\in S} v_+$.  Then $T$ is an anticlique, and $S\sub T$, so $S$ is also an anticlique, and $S\in\bigcap_{v\in R}v_+^{\epsilon_v}$.  

\end{proof}

If the hypergraph is somewhat special, we have more:
\begin{theorem}\label{thm:hyperedgesize}For any hypergraph $\G=\left<V,E\right>$ where $E\sub\left[V\right]^{\leq n}$, $\BA\G$ is $n$-free.\end{theorem}
\begin{proof}We show that $V_+$ is $n$-independent.

From the previous theorem, we need only show that $\pc2_n$ holds for $V_+$.  Let $F$ be a finite subset of $V$ such that $\prod F_+=0$.  Using the observation that $\prod F_+$ is the set of anticliques that include $F$, $F$ is not an anticlique.  Thus some hyperedge $e$ is a subset of $F$.  Then $\prod e_+=0$ as no anticlique can include that hyperedge.  Since all hyperedges have at most $n$ vertices, $\left|e\right|\leq n$, which is what we wanted.  Note that if $n=max_{e\in E}\left|e\right|$, that is, there is an edge $g$ of size $n$, that $V_+$ is not $\pa{n-1}$-independent; $\prod g =0$, but $g$ has no subset of size $n-1$ with product $0$.
\end{proof}

We have defined $\BA\G$ in two ways when $\G$ is a graph.  The two definitions are formally different; we originally used characteristic functions on $G$ rather than subsets of $G$.  However, the algebras produced are isomorphic, so there should not be any reason for confusion.

We also reverse this construction, again with the previous definition now a special case for $n=2$. Given a \ba{} $A$ with an $\omega$-independent generating set $H$, we construct a hypergraph $\G$ such that $A\cong\BA\G$; we call it the $\perp$-hypergraph of $A,H$.  The vertex set is $H$, and the hyperedge set is defined as follows; a subset $e$ of $H$ is a hyperedge if and only if the following three conditions are all true:
\begin{enumerate}
\item $e$ is finite.
\item $\prod e=0$.
\item If $f\subsetneq e$, then $\prod f\neq 0$.
\end{enumerate}

We have only finite hyperedges in this graph, and no hyperedge is contained in another.  Note that if $H$ is $n$-independent, the hyperedge set is included in $\left[H\right]^{\leq n}$.

\begin{theorem}\label{thm:hypergraphofalg}
Let $n$ be a positive integer or $\omega$, $X\sub A$ be $n$-independent and generate $A$, and $\G=\left<X,E\right>$ be the
$\perp$-hypergraph of $A$.  Then $A\cong\BA\G$.\end{theorem}
\begin{proof}Let $f:X\rightarrow X_+$ be defined so that $v\mapsto v_+$ for $v\in X$.  We claim that $f$ is an
$n$-preserving function.  If $G\sub X$ is of size $\leq n$ such that
$\prod G =0$, then it has a subset $G'$ minimal for the property of having $0$
product; thus $G'\in E$, so that $\prod G'_+=0$, and so $\prod
f\left[G\right]=0$.

$f$ is bijective, and its inverse is also $n$-preserving; the image of
$f$ is a generating set, so that $f$ extends to an isomorphism.
\end{proof}

\begin{definition}Let $\G_i\left<V_i,E_i\right>$ be hypergraphs for $i\in\s{0,1}$.  A hypergraph homomorphism is a function $f:V_0\rightarrow V_1$ such that if $e\in E_0$, then $f\left[e\right]\in E_1$.
\end{definition}

Notice that a graph homomorphism is a hypergraph homomorphism when the graphs are considered as 2-uniform hypergraphs.

With this definition of homomorphism, we have an equivalence of categories.  Let $\mathfrak{HG}_0$ be the subcategory of the category of hypergraphs including exactly those hypergraphs $\G=\left<V,E\right>$ such that:
\begin{enumerate}
\item All $e\in E$ are finite.
\item If $e\in E$ and $f\subsetneq e$, then $f\notin E$.
\end{enumerate}

We show that the category of $\omega$-free \bas{} with distinguished $\omega$-independent generating set (That is, pairs $\pa{B,G}$) is equivalent to the category $\mathfrak{HG}_0$. 
The proof is much the same as that for graphs and $2$-free \bas{}.  We let
$\mathfrak K\pa{\pa{B,G}}$ be the hypergraph on $G$ where a subset $e$ of $G$ is a
hyperedge if its $\prod e=0$ and $\prod e'\neq 0$ if $e'\subsetneq e$.  Then
$\BA{\mathfrak K\pa{\pa{B,G}}}$ is isomorphic to $\pa{B,G}$ and so on.

\section{Hypergraph Spaces}

The dual spaces to $\omega$-free \bas{} are also interesting.  Like with
graphs, a hypergraph space may be defined in terms of a hypergraph--the
definition generalizes that of a graph space.

\begin{definition}
Let $\G=\left<G,E\right>$ be a hypergraph and  $A\pa\G$ its set of anticliques.  For each $v\in
G$, we define $v_+\defeq\s{A\in A\pa\G : v\in A}$ and $v_-\defeq\s{A\in A\pa\G
: v\notin A}$.

Then the hypergraph space of $\G$ is the topology on $A\pa\G$ with $\bigcup_{v\in G} \s{v_+,v_-}$ as a closed subbase.

Any topological space $\T$ for which there is a hypergraph $\G$ such that $\T$ is homeomorphic to the hypergraph space of $\G$ is called a hypergraph space.
\end{definition}

\begin{theorem} The Stone dual of an $\omega$-free \ba{} is a hypergraph space.
\end{theorem}
\begin{proof}  Let $A$ be an $\omega$-free \ba{}.  Thus by theorem \ref{thm:hypergraphofalg}, there is a hypergraph $\G$ such that $A\cong \BA\G$.  Let $\T$ be the hypergraph space of $\G$.  We claim that $\Clop\pa\T\cong A$.

In fact, $\Clop\pa\T = \BA\G$.  On both sides here, elements are sets of anticliques of $\G$.  As $\T$ is defined by a clopen subbase, elements of $\Clop\pa\T$ are finite unions of finite intersections of elements of that subbase $\bigcup_{v\in G} \s{v_+,v_-}$.  Elements of the right hand side are sums of elementary products of elements of $\bigcup_{v\in G} \s v_+$, that is, sums of finite products of elements of $\bigcup_{v\in G} \s{v_+,v_-}$.  As the operations are the usual set-theoretic ones on both sides, they are in fact the same algebra.  

The topological result follows by duality.
\end{proof}

We repeat a few definitions from Bell and van Mill \cite{BellvanMill80} needed for some topological applications.

\begin{definition}Let $n\in \omega$ for all these definitions.

A set $S$ is $n$-linked if every $X\in\left[S\right]^n$ has
non-empty intersection.

A set $P$ is $n$-ary if every $n$-linked subset of $P$ has non-empty
intersection.

A compact topological space $\T$ has compactness number at most $n$, written $\cmpn\pa\T \leq n$,
if and only if it has an $n$-ary closed subbase.
$\T$ has compactness number $n$, written $\cmpn\pa\T = n$, if and only if $n$ is the least integer
for which $\cmpn \pa \T \leq n$. $\cmpn\pa\T = \omega$ if there is no such $n$.
\end{definition}

The following generalizes and algebraizes proposition 3.1 of Bell \cite{Bell82}.

\begin{proposition}\label{prop:bell}If a \ba{} $A$ is $n$-free for some $2\leq n \leq \omega$, then $\cmpn\pa{\Ult A}\leq n$.\end{proposition}
\begin{proof} This is vacuous if $n=\omega$.  If $n<\omega$, then $\Ult\pa A$ is a hypergraph space for a hypergraph $\G$ with all hyperedges of size $\leq n$.

We take the clopen subbase $S=\bigcup_{v\in G} \s{v_+,v_-}$ of the hypergraph space of $\G$ and show that it is $n$-ary.  Let $\mathcal F\sub S$ be $n$-linked.  We may write $\mathcal F=\s{v_+ : v\in A}\cup \s{v_- : v\in B}$ for some $A,B\sub G$.  Since $v_+\cap v_-=\emptyset$ and $n\geq 2$, $A\cap B=\emptyset$.  Let $A'$ be a finite subset of $A$.  Since any product of $n$ or fewer elements of $\mathcal F$ is non-zero, $A'$ must be an anticlique in $\G$; if not, then $\prod A'_+=0$, so then $A'_+$ would have a subset of size $n$ with empty intersection, contradicting that $\mathcal F$ is $n$-linked.
Thus $\bigcap \s{v_+: v\in A'}\in \bigcap \mathcal F$, that is, $\mathcal F$ has non-empty intersection and thus $S$ is $n$-ary.
\end{proof}

Bell's \cite{MR803933} corollary 5.2 shows that certain topologies on
$\left[\omega_1\right]^{\leq m}$ have compactness number $n$ for certain 
$n,m\leq \omega$.  These topologies are the hypergraph spaces of
$\left<\omega_1,\left[\omega_1\right]^{2n-3}\right>$ and
$\left<\omega_1,\left[\omega_1\right]^{2n-2}\right>$.
\begin{theorem}For infinitely many $n\in\omega$, there is a \ba{} which is $n$-free and is not $\pa{n-1}$-free.
\end{theorem}
\begin{proof}Let $k$ be the least integer for which $\BA{\left<\omega_1,\left[\omega_1\right]^{2n-3}\right>}$ is $k$-free and $\ell$ be the least integer for which $\BA{\left<\omega_1,\left[\omega_1\right]^{2n-2}\right>}$ is $\ell$-free. We have that
$n\leq k \leq 2n-3$ and $n\leq \ell\leq 2n-2$.  The lower bounds are a
consequence of the compactness numbers of those spaces (Bell's \cite{MR803933} result and proposition \ref{prop:bell}), while the upper bounds
are a consequence of theorem \ref{thm:hyperedgesize}.

Thus we have, for arbitrary $n\in\omega$, an $\omega$-free \ba{} of finite freeness at least
$n$.
\end{proof}

\section{Algebraic Constructions}\label{sec:ac}

\newcommand{\GS}{$\mathscr{GS}$ }
In this section, we consider the categories of $n$-independently generated \bas{} and of hypergraph spaces with respect to closure under algebraic constructions.

Closure of a class of \bas{} under homomorphic images is equivalent to closure of
their dual spaces under closed subspaces.  Every \ba{}  is the homomorphic image of a free
one, thus the homomorphic image of an $\omega$-free one.  Since there are
non-$n$-free \bas{} for every $n$, the class of $n$-free \bas{} is not closed under
homomorphic images, and the category of graph spaces \GS is not closed under closed subspaces.
 
We will show in section \ref{sec:cf} that complete \bas{} are not $\omega$-free.
As $\P\left(\kappa\right)$ is isomorphic to $\leftexp \kappa 2$, the class of $\omega$-free \bas{} is not closed under infinite products, and 
\GS is not closed under infinite free products (i.e. Stone-\v{C}ech compactifications of disjoint unions). 
However, \PFBA is closed under some product constructions.

\begin{theorem}\label{finprod}Let $2\leq n\leq \omega$. If $H\sub A$ and $K\sub B$ are $n$-independent, then  
$L\defeq\left(H\times \left\{0\right\}\right) \cup \left(\left\{0\right\}\times K\right)$ is $n$-independent in $A\times B$.
\end{theorem}
\begin{proof}

We will apply proposition \ref{prop:a}.  Suppose that $F\in\subs H {<\omega}, G\in\subs K {<\omega}, \epsilon\in \leftexp F2, \delta\in\leftexp G2$, and $\prod_{x\in F}\pa{x,0}^{\epsilon_x}\cdot\prod_{y\in G}\pa{0,y}^{\delta_y}=0$.  If there are $x\in F$ and $y\in G$ such that $\epsilon_x=\delta_y=1$, then $\pa{x,0}\cdot\pa{0,y}=0$ as desired.  Otherwise, without loss of generality, we may assume that $\epsilon\left[F\right]\sub\s0$.  Then $\prod_{x\in F}\pa{x,0}^{\epsilon_x}=\pa{\prod_{x\in F}-x,1}$, so that $\prod_{y\in Y}y^{\delta_y}=0$; then the $n$-independence of $K$ gives the result.
\end{proof}

It is important to note that $L$ does not generate $\left<H\right>\times \left<K\right>$; in fact (theorem \ref{thm:fcf_not2free}), the product of $n$-free \bas{} is not in general $n$-free. However, it is
the case that $\left<H\right>\times \left<K\right>$ is a simple extension of the subalgebra generated by $L$; $\left<L\right>\pa{\pa{1,0}}=\left<H\right>\times \left<K\right>$.

This result generalizes to infinite products quite easily, though the notation
is considerably more cumbersome.

\begin{theorem}\label{bigbindset} For $2\leq n \leq \omega$, if $\left<A_i : i\in I\right>$ is a system of \bas{} and for every
$i\in I$, $H_i\sub
A_i$ is $n$-independent in $A_i$, then the set $H\defeq\bigcup_{i\in I}
p_i\left[H_i\right]$, where
$$p_i\left(h\right)\left(j\right)\defeq\left\{\begin{array}{ll}h & i=j \\ 0 & i\neq
j\end{array}\right.$$ is $n$-independent in $A\defeq\prod_{i\in I}A_i$ and
$\prod^{\mathrm w}_{i\in I}A_i$.
\end{theorem}
\begin{proof}

This is essentially the same as  theorem \ref{finprod} with
more cumbersome notation.  

$p_i\left(h\right)$ is the function in
$A$ that is $0$ in all but the $i$th coordinate and is $h$ in the $i$th
coordinate, so that the projections $\pi_i\left[p_i\left[H_i\right]\right]=H_i$ and for $i\neq
j$,  $\pi_j\left[p_i\left[H_i\right]\right]=\left\{0\right\}$.

We apply proposition \ref{prop:a}.  Suppose that $R\in \subs H {<\omega}, \epsilon\in\leftexp R2$, and $\prod_{x\in R}x^{\epsilon_x}=0$.  Let $J\defeq\s{i\in I:R\cap p_i\left[H_i\right]\neq \emptyset}$.  If $J$ is a singleton, say $J=\s i$, then the $n$-independence of $H_i$ clearly makes $H$ $n$-independent. So we now concern ourselves with the case that  $\left|J\right|>1$, that is, we have distinct $i,j\in J$. If there are $x\in p_i\left[H_i\right]$ and $y\in p_j\br{H_j}$ with $\epsilon_x=\epsilon_y=1$, then $x\cdot y=0$ and we have our conclusion.  So we may assume that there is at most one $i\in J$ for which there is an $x\in p_i\br{H_i}$ such that $\epsilon_x=1$.  Then $\prod\s{x^{\epsilon_x}:x\in R, x\notin p_i\br{H_i}}$ has $i$-th coordinate $1$, and so
$$0=\prod_{x\in R}x^{\epsilon_x} = \prod\s{x^{\epsilon_x}:x\in R\cap p_i\br{H_i}},$$
 and the $n$-independence of $H_i$ make $H$ $n$-independent.

\end{proof}

When $n=2$, we can also consider the $\perp$-graph and intersection graph of $H$ in the above theorem.
The intersection graph is easily described: two elements of $H$ have non-zero
product if and only if they have non-zero product in one of the factors, so that the
intersection graph is the disjoint union of the intersection graphs of the
$H_i$.  The $\perp$-graph is more complex.  The $\perp$-graph of each $H_i$ is
an induced subgraph, but these subgraphs are connected to each other--each
vertex in $H_i$ is connected to every vertex in $H_j$ for $i\neq j$. This
construction is the ``join".

In other words, for any collection $\left<\G_i\right>$ of graphs,  $\BC{\bigcup_{i\in I}\G_i}\leq \prod_{i\in I}\BC{\G_i}$ and 
$\BA{\biguplus_{i\in I}\G_i}\leq \prod_{i\in I}\BA{\G_i}$.

The use of the word ``free" in $n$-free is warranted by the following:

\begin{theorem}\label{freeprod}  Suppose that $A\defeq\bigoplus\limits_{\stackrel{C}{i\in I}} A_i$ is an amalgamated free product of subalgebras $A_i$ for $i\in I$, where $C\leq A_i$ for each $i\in I$, $A_i\cap A_j =C$ for $i\neq j$, $A_i$ is $n$-free over $H_i$, and $C\leq\left<H_i\cap H_j\right>$.  Then $A$ is $n$-free over $\bigcup_{i\in I}H_i$.

\end{theorem}
\begin{proof}For convenience, assume that each $A_i \leq A$, $C\leq A_i$ and that, for
$i\neq j$, $A_i\cap A_j = C$, and
that $H_i$ is a set over which $A_i$ is $n$-free.  
We show that $A$ is $n$-free over $H\defeq\bigcup_{i\in I} H_i$.  

Let $B$ be a \ba{}, and  $f:H\rightarrow B$ be $n$-preserving.  Then
for each $i\in I$, $f_i:=f\res H_i$ is also $n$-preserving.  So each $f_i$ extends to a unique homomorphism $\phi_i : A_i\rightarrow B$.  That $\phi_i\res C = \phi_j\res C$ is clear as $C\sub\left<H_i\cap H_j\right>$.

 Then the universal property of amalgamated free products gives a unique homomorphism $\phi : A\rightarrow B$ that extends every $\phi_i$.  Note that
$$\phi\res H = \phi \res\bigcup_{i\in I} H_i=\bigcup_{i\in I}\left(\phi\res H_i\right)=\bigcup_{i\in I} \left(\phi_i\res H_i\right) = \bigcup_{i\in I} f_i = f.$$
So we have a unique extension of $f$ to a homomorphism, which is what we
wanted.\end{proof}

This of course includes free products, so
dually, we have that \GS is closed under arbitrary products.

An example where $C\neq\s{0,1}$ is as follows:
Let $\G$ be the complete graph on the ordinal $\omega_1 +\omega$ and $\H$ the complete graph on the ordinal interval $\left(\omega_1,\omega_1\cdot 2\right)$.  Then $\BA\G\cong\BA\H\cong\FR{\omega_1}$.  Note that $G\cap H=\left(\omega_1,\omega_1+\omega\right)$ so that $G_+\cap H_+=\left(\omega_1,\omega_1+\omega\right)_+$; we let $C=\left<\left(\omega_1,\omega_1+\omega\right)_+\right>\cong\FR\omega$.  It is clear that $C$ is as required in theorem \ref{freeprod}. Then we have that $\BA\G\mathop{\oplus}\limits_C \BA\H$ is $2$-free over $G_+\cup H_+$.

If $C$ is $2$-free over $\bigcap_{i\in I}H_i$, the $\perp$-graph of $\bigcup_{i\in I}H_i$ is easily described in terms of those of
$H_i$. Indeed, our equivalence of categories implicitly does this already--it
is the ``amalgamated free product'' or  ``amalgamated disjoint union'' in the category of graphs--i.e. the same universal property holds.  More concretely, given a set of graphs $\G_i=\left<G_i,E_i\right>$, each of which has $\mathscr F=\left<F,E\right>$ as a  subgraph, the amalgamated disjoint union of the $\G_i$ over $\mathscr F$ is a graph on the union of the vertex sets where two vertices are adjacent if and only if they are adjacent in some $\G_i$.  That is, elements of $G_i\sm F$ and $G_j\sm F$ are not adjacent for $i\neq j$.  

In case $C=2$ and we have a free product, the $A_i$ form a family of independent subalgebras, so two elements
of $H$ (constructed in the proof above) have product zero if and only if they are in the same $H_i$ and have zero product
in $A_i$.  So the $\perp$-graph of $H$ is the disjoint union of the $\perp$
graphs of the $H_i$.
The intersection graph of $H$ is similarly constructed from those of the $H_i$:
the independence of the $A_i$ means that the intersection graph of $H$ is the
join of the intersection graphs of the $H_i$.

That is, $\bigoplus_{i\in I}\BA{\G_i} = \BA{\bigcup_{i\in I} \G_i}$ and 
$\bigoplus_{i\in I}\BC{\G_i} = \BC{\biguplus_{i\in I}\G_i}$.

We now give the examples promised at the end of section \ref{sec:gs}.  In that section, we topologically showed that $\BA{K_4}\cong\BA{P_3}\cong\s{0,1}^5$, where $K_4$ is the complete graph on 4 vertices and $P_3$ is the path on three vertices.  So if we let $\kappa$ be an infinite cardinal, $\bigcup_{\alpha<\kappa}K_4 \not\cong \bigcup_{\alpha<\kappa}P_3$, but $\BA{\bigcup_{\alpha<\kappa}K_4}\cong\BA{\bigcup_{\alpha<\kappa}P_3}$, as both are isomorphic to $\bigoplus_{\alpha<\kappa}\s{0,1}^5$.

Products of $n$-free \bas{} behave in a somewhat more complicated manner.  As discussed previously, infinite products
of $\omega$-free \bas{} are not necessarily $\omega$-free.  

\begin{theorem}\label{thm:fcf_not2free}$\FinCo{\omega_1}\times\FR{\omega_1}$ is not $2$-free.\end{theorem} 
\begin{proof} We use subscript function notation for the coordinates of tuples; i.e. $\pa{a,b}_0=a$ and $\pa{a,b}_1=b$.  We also extend this to sets of tuples; $\s{\pa{a,b},\pa{c,d}}_0=\s{a,c}$.

  We proceed by contradiction; suppose that $A\defeq\FinCo{\omega_1}\times\FR{\omega_1}$ is $2$-free over $X$, where $0\notin X$, that is, $X$ is $2$-independent. 

Consider $a_\alpha\defeq\pa{\s\alpha,0}$ for $\alpha<\omega_1$.  $a_\alpha$ is
an atom in $A$, so it must be an elementary product of $X$, that is,
$a_\alpha=\prod_{x\in H_\alpha}x^{\epsilon\pa{\alpha,x}}$, with
$H_\alpha\in\left[X\right]^{<\omega}$.  So let
$M\in\left[\omega_1\right]^{\omega_1}$ be such that $\s{H_\alpha:\alpha\in M}$
is a $\Delta$-system with root $F$.  Let $G_\alpha\defeq H_\alpha\sm F$.
Since $M=\bigcup_{\delta\in \leftexp F2}\s{\alpha\in M:
\forall x\in F\left[\epsilon\pa{\alpha,x}=\delta_x\right]}$, there is an uncountable $N\sub M$
such that $\epsilon\pa{\alpha,x}=\epsilon\pa{\beta,x}$ for all $\alpha,\beta\in
N$ and all $x\in F$, so that we may write, for $\alpha\in N$, $a_\alpha=\prod_{x\in
F}x^{\delta_x}\cdot \prod_{x\in G_\alpha}x^{\epsilon\pa{\alpha,x}}$. 
For each $\alpha\in N$, let $G'_\alpha\defeq\s{x\in G_\alpha :
\epsilon\pa{\alpha,x}=1}$.  If $\alpha,\beta\in N$ with $\alpha\neq\beta$, then
there are $x\in G'_\alpha$ and $y\in G'_\beta$ such that $x\cdot y=0$, by the
$2$-independence of $X$ and the fact that $0=a_\alpha\cdot a_\beta=\prod_{x\in F}x^{\delta_x} \cdot\prod_{x\in G_\alpha}x^{\epsilon\pa{\alpha,x}}\cdot\prod_{x\in G_\beta}x^{\epsilon\pa{\beta,x}}$,  thus $\prod G'_\alpha \cdot \prod G'_\beta=0$.   Since $\FR{\omega_1}$ has cellularity $\omega$, the set $\s{\alpha\in N:\pa{\prod G'_\alpha}_1\neq 0}$ is countable, hence $P\defeq N\sm \s{\alpha\in N:\pa{\prod G'_\alpha}_1\neq 0}$ is uncountable and for $\alpha\in P$, $\pa{\prod G'_\alpha}_1=0$.  Since $\pa{\prod G'_\alpha}_0\cdot\pa{\prod G'_\beta}_0=0$ for distinct $\alpha,\beta\in P$, each $\pa{\prod G'_\alpha}_0$ is finite when $\alpha\in P$.

$X$ must generate $\pa{1,0}$; let $b_j$ for $j<n$ be disjoint elementary
products of $X$ such that $\sum_{j<n}b_j=\pa{1,0}$.  Thus there must be exactly
one $i<n$ such that $b_{i0}$ is cofinite; without loss of generality, $i=0$ so that $b_{00}$ is cofinite
and $b_{01}=0$.  Write $b_0$ as an elementary product, that is
$b_0=\prod_{j<n}c_j^{\xi_j}$ with each $c_j\in X$.  Then choose an $\alpha\in
P$ such that $\prod G'_\alpha\leq b_0$ and $G'_\alpha\cap\s{c_j:j<n}=\emptyset$.
Then $\prod G'_\alpha\cdot\sum_{j<n}c_j^{1-\xi_j}=0$, so $\rng\xi=\s0$; that
is, $b_0=\prod_{j<n} -c_j$

Note that $X_1$ generates $\FR{\omega_1}$, so it must be uncountable, thus\\  
$\pa{X\sm\s{c_j:j<n}\sm F}_1$ is also uncountable; let $Y\sub X$ be such that $Y_1$ is
an uncountable independent subset of $\pa{X\sm\s{c_j:j<n}\sm F}_1;$ such a $Y$
exists by theorem 9.16 of Koppelberg \cite{Handbook}.  Note that no finite
product of elements of $Y$ is $0$. Let
$\theta:Y\rightarrow\s{0,1}$ be such that $d_y\defeq\pa{y^{\theta_y}}_0$ is finite for
each $y\in Y$.

Consider $\s{d_y:y\in Y}$; Each $d_y$ is finite and $Y$ is an uncountable set, and thus there is 
an uncountable $Z\sub Y$ where $\s{d_y:y\in Z}$ is a $\Delta$-system with root $r$.  Let $y,z,t\in Z$ be 
distinct.  Then let $e_y\defeq d_y\sm r$, $e_z\defeq d_z\sm r$, and $e_y\defeq d_z\sm r$.  Then $d_y\cdot d_z \cdot -d_t=\pa{e_y \cup r}\cap \pa{e_z\cup r}\cap\pa{\omega_1\sm\pa{e_t\cup r}}=r\cap \pa{\omega_1\sm\pa{e_t\cup r}}=\emptyset$.  Then $\prod_{j<n}-c_j\cdot y^{\theta_y}\cdot z^{\theta_z}\cdot t^{1-\theta_t}=0$ and again, the only elements with exponent $1$ are elements of $Y$ and thus there is no disjoint pair, contradicting proposition \ref{prop:a}.

So we have a contradiction and thus there is no $2$-independent generating set
for $A$.
\end{proof}
 This is also an example of a simple extension of a
$2$-free \ba{} that is not $2$-free; the full product is a simple extension by
$\left(0,1\right)$ of the subalgebra generated by the set in theorem
\ref{finprod}.

The dual of this theorem is that we have two graph spaces whose disjoint union
is not a graph space; in fact we can say a bit more since the disjoint union of two
supercompact spaces is supercompact.  We show a slightly more general result
here:

\begin{proposition}\label{prop:union-scomp}If $X$ and $Y$ are $n$-compact
spaces, then $X\dot\cup Y$ is $n$-compact.\end{proposition}
\begin{proof}Suppose that $S$ and $T$ are $n$-ary subbases for the closed sets of $X$ and $Y$
respectively; that is, for any $S'\sub S$ with $\bigcap S'=\emptyset$, there
are $n$ members $a_1,a_2,\ldots,a_n$ of $S'$ such that $a_1\cap
a_2\cap\ldots\cap a_n=\emptyset$, and similarly for $T$.  Then $W\defeq S\cup
T\cup\s{X,Y}$ is an $n$-ary subbase for the closed sets of $X\dot\cup Y$.\end{proof}

So, letting $n=2$, the dual space of $\FinCo{\omega_1}\times\FR{\omega_1}$ is supercompact, but
is not a graph space.

\begin{theorem}$\FR{\omega_1}\times\FinCo{\omega_1}$ is $3$-free.\end{theorem}
\begin{proof} Let $\s{x_\alpha:\alpha<\omega_1}$ be an independent generating
set for $\FR{\omega_1}$.  Then the set
$X\defeq\s{\pa{x_\alpha,\s\alpha}:\alpha<\omega_1}\cup\s{\pa{0,1}}$ is a
$3$-independent generating set for $\FR{\omega_1}\times\FinCo{\omega_1}$.
That $X$ generates $\FR{\omega_1}\times\FinCo{\omega_1}$ is clear.  We use
proposition \ref{prop:a} to show that $X$ is $3$-independent.  Take any
$R\in\left[X\right]^{<\omega}$ and $\epsilon\in \leftexp R2$ such that $\prod_{x\in R}x^{\epsilon_x}=\pa{0,0}$.  Since there is no elementary product of elements of $\s{x_\alpha:\alpha<\omega_1}$ that is $0$, $\pa{0,1}\in R$ and $\epsilon_{\pa{0,1}}=1$.  Then there is a pair $a,b$ of elements in $R$ such that $\pi_2\pa a \perp \pi_2 b$ and $\epsilon_a=\epsilon_b=1$, so that $\s{\pa{0,1},a,b}\sub R$ and $\pa{0,1}\cdot a\cdot b =0$.  
\end{proof}

\section{Cardinal Function Results}\label{sec:cf}
Cellularity and independence have been considered earlier.  Here we give a few
results relating other cardinal functions to properties of $\perp$-graphs and
intersection graphs.  We will always assume that the graphs and algebras are infinite in
this section.

It is worth noting that the normal form of an element of a $2$-free \ba{}  (finite sums of finite elementary products) can be translated somewhat succinctly into graph language.  We give these in terms of clique algebras, but adding ``anti'' to ``clique'' everywhere is all that is necessary to change to anticlique algebras.  An elementary product has two finite sets of generators associated with it, those that occur in it without being complemented, and those that occur complemented.  As vertices, this means that an elementary product is the set of all cliques that include the clique of all vertices of the first kind (if these vertices are not a clique in their own right, their product is empty) and omits all the vertices of the second kind. Then any element of the algebra is a finite union of such sets.

The following is useful in several places:

\begin{lemma}\label{finco_is_image}Let $A$ be $\omega$-free over $G$ and $\omega
\leq\kappa=\left|A\right|$.  Then $B\defeq\FinCo\kappa$ is a homomorphic image of $A$.
\end{lemma}
\begin{proof}
Any bijective function $f:G\to \At B$ is $\omega$-preserving as all elements of $\At B$ are disjoint.  Since $A$ is $\omega$-free, $f$ extends to a homomorphism $\tilde f$ from $A$ to $B$.  Since the image of $f$ includes a set of generators, $\tilde f$ is surjective as well; that is, $B$ is a homomorphic image of $A$.
\end{proof}

The first use of this is that no infinite $\omega$-free \ba{} has the countable separation property. The \CSP{} is inherited by homomorphic images (5.27(c) in Koppelberg \cite{Handbook}), so if any infinite $\omega$-free \ba{} of size $\kappa$ 
has the \CSP{}, then by \ref{finco_is_image}, $\FinCo\kappa$ has the \CSP{}, which is a contradiction.  In particular, $\P\pa\omega /\mbox{fin}$ is not $\omega$-free.

\subsection{Spread}\label{spread}
We show that the spread of an $\omega$-free \ba{}  is  equal to its cardinality.

Theorem 13.1 of Monk \cite{CINV2} gives several equivalent definitions of spread, all of which have the same attainment properties; the relevant one to our purposes is the following.
$$s\pa A=\sup\left\{c\pa B : B\mbox{ is a homomorphic image of }A\right\}.$$

\begin{theorem}For $A$ $\omega$-free, $s\pa A=\left|A\right|$.  Furthermore, it is attained.\end{theorem}
\begin{proof}
From lemma \ref{finco_is_image}, $B=\FinCo{\left|A\right|}$ is a homomorphic image of $A$.  Since $c\pa B=\left|B\right|=\left|A\right|$, an element of the set in the above definition of $s\pa A$ is $\left|A\right|$.  Thus $s\pa A=\left|A\right|$ is attained.
\end{proof}

Another equivalent definition of spread is 
$$s\pa A=\sup\s{\left|X\right| : X\mbox{ is ideal-independent.}}$$
Here the definition of ideal independence is $0,1\notin X$ and $\forall x\in X,
x\not\in\left<X\sm\s x\right>^{id}$; equivalently, for distinct $x,x_1,\ldots,x_n\in X$,
$x\not\leq x_1+x_2+\ldots +x_n$. For $\omega$-free \bas{}, it is already
shown that $s\pa A$ is attained in every sense; we further show that in this sense, it is
attained by an $\omega$-independent generating set.  
\begin{theorem}\label{pindisidind} If $H$ is $\omega$-independent, then $H$ is
ideal-independent.
\end{theorem}
\begin{proof}
Let $x,x_1,x_2,\ldots,x_n\in H$ be distinct.  Let $F=\s x$ and
$G=\s{x_1,x_2,\ldots,x_n}$.  $\prod F\neq 0$ and $F\cap G=\emptyset$ so that \pc3 implies that $\prod F\not\leq \sum G$, that is, $x\not\leq x_1+x_2+\ldots x_n$.  Thus $H$ is ideal-independent.
\end{proof}

As they are greater than or equal to $s$, Inc, Irr, h-cof, hL, and hd are also equal to cardinality for $\omega$-free \bas{}.
This result also determines that $\left|\mbox{Id} A\right|=2^{\left|A\right|}$
as $2^{sA}\leq \left|\mbox{Id} A\right|$.  Then since $s$ is attained, $\left|\mbox{Sub}
A\right|=2^{\left|A\right|}$ as well.

There is more to say about incomparability and irredundance; they are also attained by the $\omega$-free generating set.

\begin{proposition}\label{incomp}
If $H$ is a $\omega$-independent set, then $H$ is incomparable.
\end{proposition}
\begin{proof}Assume that $g,h\in H$ are such that $g<h$.  Let $F\defeq\s g$ and $G\defeq\s h$.  Then $0\neq \prod F \leq \sum G$ and $F\cap G=\emptyset$, contradicting $\pa{\perp3}$.  
\end{proof}

\begin{theorem} If $H$ is a $\omega$-independent set, then $H$ is irredundant.
\end{theorem}
\begin{proof}Let $B\defeq\left<H\right>$.

Assume not, namely let $a\in H$ be such that
$\left<H\sm\left\{a\right\}\right>= B$.  Then $a$ can be written as a finite
sum of elementary products of elements of $H\sm\left\{a\right\}$, say $\sum_{i\leq
n} \prod_{j\leq m_i} a_{ij}^{\epsilon_{ij}}$.  We construct an  
$f:H\rightarrow 2$ as follows. Let $f$ be identically zero on $H\sm\left\{a\right\}$.  
Define \begin{equation}\label{irreqn} f\left(a\right)=-\sum_{i\leq
n} \prod_{j\leq m_i} 0^{\epsilon_{ij}},\end{equation} that is, $1$ if every $\epsilon_{ij}$
is $0$, and $0$ otherwise.  Then $f$ does not
extend to a homomorphism as the complement of  the right hand side of equation \ref{irreqn} is $f$ applied to the finite sum of elementary products that produces $a$.  However, there is a contradiction as $f$ is
$\omega$-preserving.

\end{proof}

\subsection{Character}

The character of an $\omega$-free \ba{}  is yet one more invariant equal to cardinality.  Namely, at the bottom of page 183 in Monk \cite{CINV2}, it is shown that if $A$ is a homomorphic image of $B$, then $\chi\pa A \leq \chi \pa B$.  For $B$ $\omega$-free, let $A=\FinCo{\left|B\right|}$, so that $A$ is a homomorphic image of $B$ by lemma \ref{finco_is_image}, so we have that $\left|B\right|=\chi\pa A \leq \chi\pa B \leq \left|B\right|$.

\subsection{Length}

We claim that the length (and therefore depth) of an $\omega$-free \ba{} is $\aleph_0$.  This
uses several preceding results.

\begin{theorem} If $A$ is $\omega$-free, then $A$ has no uncountable chain.
\end{theorem}
\begin{proof}Let $A$ be $\omega$-free over $G$.

Recall from theorem \ref{semigroup} that $A$ is a semigroup algebra over the set $H$ of finite products of elements of $G\cup\left\{0,1\right\}$.
 For $h\in H\sm\left\{0,1\right\}$, choose $g_1,\ldots,g_n\in G$ such that
$h=g_1\cdot\ldots\cdot g_n$ and set $h_G\defeq\left\{g_1,\ldots,g_n\right\}.$

Due to the result of Heindorf \cite{MR1155389}, if there is an uncountable chain in $A$,
there is an uncountable chain in $H$.  So by way of contradiction, we assume that there is an
uncountable chain $C\sub H$.  Without loss of generality, we may assume that
$0,1\not\in C$ so that every element of $C$ is a finite product of elements of
$G$. 

Let $C_G\defeq\left\{h_G:h\in C\right\}$.  We note that $$\bigcup C_G=\bigcup_{h\in
C}h_G\sub G$$ is the set of all elements of $G$ that are needed to generate the
elements of $C$, that is,

$C\sub \left<\bigcup C_G\right>$. so $C$ is a chain in that subalgebra of $A$
as well.

In order to reach a contradiction, we first show
that there are no finite subsets of $\bigcup
C_G$ with zero product.
Take $F\in\left[\bigcup C_G\right]^{<\omega}$.
Then for each $v\in F$, there is a $c_v\in C_G$ such that $v\in c_v$.  Note
that $\prod c_v\in C$ and $\prod c_v\leq v$.
Thus $\s{\prod c_v : v\in F}\sub C$, so $0\neq\prod\s{\prod c_v : v\in
F}\leq\prod F$, and hence $\prod F\neq 0$.

Thus $\bigcup C_G$ has no finite subset with zero product.  As $\bigcup C_G\sub G$ is
$\omega$-independ-ent, by lemma \ref{pindtoind}, it is independent.  Thus $\left<\bigcup
C_G\right>$ is free and hence has no uncountable chain, contradicting our original
assumption.
\end{proof}

\subsection{Cellularity}
Results on cellularity of $\omega$-free \bas{} are necessarily somewhat
different from the preceding results.  In the previous results, cardinal functions on
$\BC\G$ and $\BA\G$ were computed using nothing more than the size of the vertex set of
$\G$.  As $\FinCo\kappa=\BA{K_\kappa}$ and $\FR\kappa=\BC{K_\kappa}$, yet have cellularity $\kappa$ and $\omega$ respectively, more
information about $\G$ is necessary--that is, the edge set matters.  From proposition \ref{ind} and the
ensuing comments, we have that
$c\pa{\BC\G} \geq \sup\s{\left|A\right|: A\;\mbox{is an anticlique in}\;
\G}+\aleph_0$ and \\$c\pa{\BA\G} \geq \sup\s{\left|C\right|: C\;\mbox{is a clique in}\; \G}+\aleph_0$.

\subsection{Independence}
Again, the edges of $\G$ must be considered in calculating
$\mbox{Ind}\pa{\BC\G}$.  

We have that
$$\mbox{Ind}\pa{\BA\G}\geq\sup\s{\left|A\right|:A\;\mbox{is an anticlique
in}\;\G}+\aleph_0\;\mbox{and}$$ $$\mbox{Ind}\pa{\BC\G}\geq\sup\s{\left|C\right|:
C\;\mbox{is a clique in}\;\G}+\aleph_0.$$

\subsection{Number of Endomorphisms}
Again, we have a ``largest possible'' result.

\begin{theorem}\label{non_atom_endo}Let $A$ be $\omega$-free over $H$. Then $\left|\mathrm{End}\; A\right|= 2^{\left|A\right|}$.
\end{theorem}
\begin{proof}
For each $x\in H$, choose $y_x\in A$ such that $y_x<x$. For each $J\subset H$, define $f_J:H\rightarrow A$ as 
$$f_J\left(x\right)=\left\{\begin{array}{ll} y_x & x\in J \\ x & \mbox{otherwise.}\end{array}\right.$$
$f_J$ is $2$-preserving and extends to an endomorphism.  So we have exhibited $2^{\left|A\right|}$ endomorphisms.
\end{proof}

\subsection{Number of Automorphisms}

\begin{theorem}\label{atom_auto}Let $A$ be $\omega$-free over $H$.
If $H$ includes  $\kappa$ atoms of $A$, then $\left|\mathrm{Aut}\; A\right|\geq 2^\kappa$.
\end{theorem} 
\begin{proof}
From proposition \ref{incomp}, $H$ is an incomparable set, so if $a\in H$ is an atom, then for any $b\in H\sm\s a$, $a\perp b$. For any permutation $\sigma$ of the atoms in $H$, let $f_\sigma :H \rightarrow H$ be such that $$f_\sigma\left(x\right)=\left\{\begin{array}{ll}\sigma\left(x\right) & x\mbox{ is an atom}\\ x & \mbox{otherwise.}\end{array}\right.$$
Note that $f_\sigma$ is a permutation of $H$ and that its inverse is $f_{\sigma^{-1}}$.  $f_\sigma$ is $2$-preserving: if $a,b$ are not atoms, then $f_\sigma$ does nothing to them, so $a\perp b$ implies $f_\sigma\left(a\right) \perp f_\sigma\left(b\right)$, whereas if $a$ is an atom, then $a\perp b$ for any $b\in H\sm\s a$. As $f_\sigma\left(a\right)$ is also an atom, and is not equal to $f_\sigma\pa b$, $f_\sigma\left(a\right)\perp f_\sigma\left(b\right)$.  We also note that this argument for $2$-preservation applies to $f_{\sigma^{-1}}=f_\sigma^{-1}$, so that we have the stronger statement that for $a,b\in H$, $a\perp b \iff f_\sigma\pa a \perp f_\sigma\pa b$.

Every $f_\sigma$ extends to a homomorphism, call it $F_\sigma$.  Then
$F_\sigma\circ F_{\sigma^{-1}}$ extends the identity on $H$, and hence is the
identity on $A$.  Similarly, $F_{\sigma^{-1}}\circ F_\sigma$ is the identity on
$A$.  So $F_\sigma$ is an automorphism.

Since there are $2^\kappa$ such $\sigma$, we have exhibited $2^\kappa$
automorphisms.

\end{proof}

Note that an atom $v_+$ in $G_+
\sub \BC\G$ corresponds to an isolated vertex. In other words, these automorphisms of $\BC\G$
correspond to automorphisms of $\G$ that do nothing but permute its isolated vertices, so this is not really that interesting.
To use the equivalence of categories, we should compare to $\BA\G$, in which case we are permuting vertices that have every possible edge.

\begin{figure}[!ht]
\label{pff}
\input{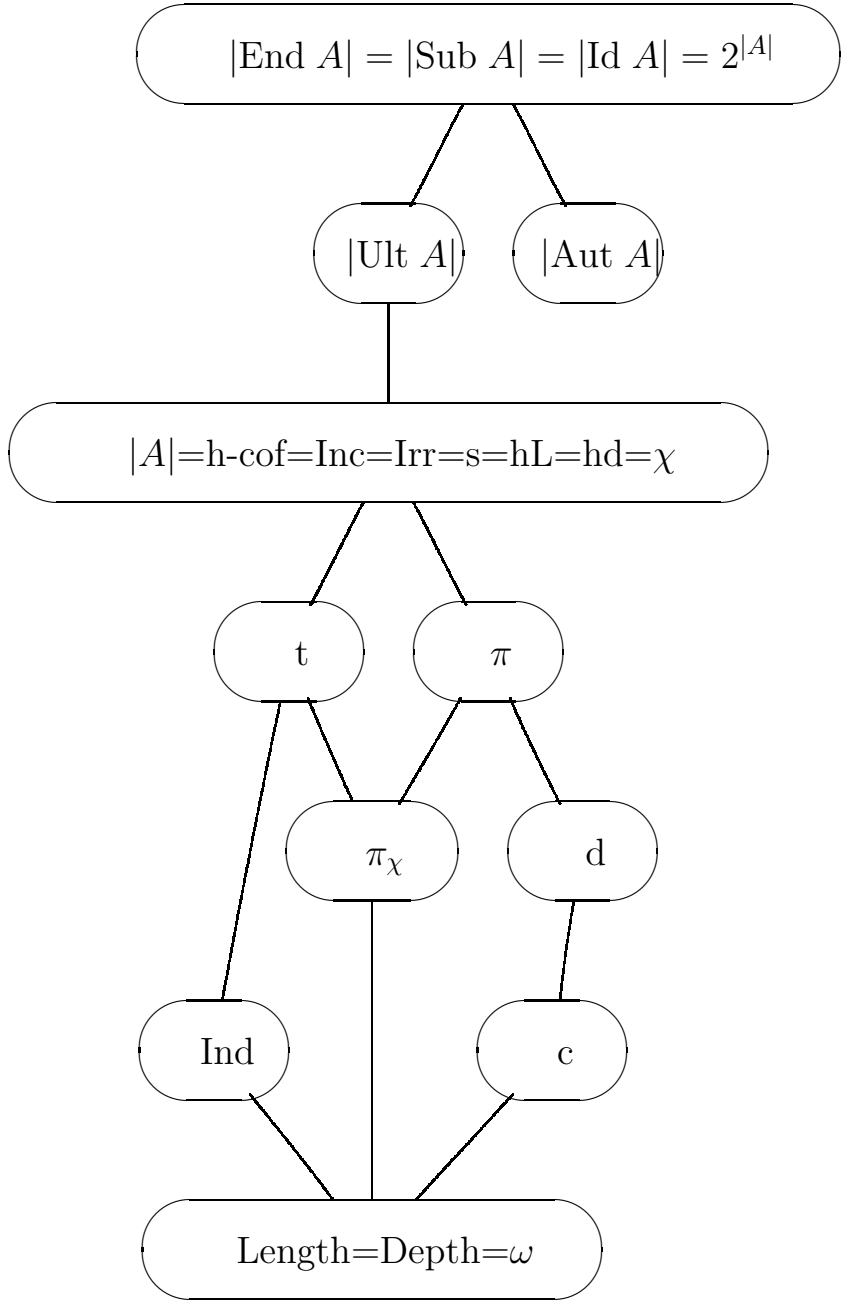}
\caption{Relationship of Cardinal Functions on $\omega$-free \bas{}.}
\end{figure}
\clearpage
Most of the lines in figure 1 indicate the possibility of a large difference.  We give examples going roughly top to bottom.

We do not know if $\left|\mathrm{Aut}\; A\right|$ may be small.

$\FinCo\kappa$ has $\kappa$ ultrafilters, so gives an example of as large as possible difference between $|\mathrm{Ult}\; A|$ and $2^{|A|}$.

$\FR\kappa$ has $2^\kappa$ ultrafilters, so gives an example of large difference between $|A|$ and $|\mathrm{Ult}\; A|$.

$\FinCo\kappa$ has tightness $\omega$, so this is a large difference between $t$ and $|A|$.

The difference between $\pi\pa A$ and $\left|A\right|$ is small.

$\FR\kappa$ has tightness $\kappa$, so gives a large difference between depth and $t$.

$\FR\kappa$ has cellularity $\omega$ and density $\lambda$ where $\lambda$ is the smallest cardinal such that $\kappa\leq 2^\lambda$, giving a large difference.

$\FinCo\kappa$ has cellularity $\kappa$ and depth $\omega$, giving a large difference there.

$\FR\kappa$ has depth $\omega$ and $\kappa$ independence, so that difference
may be large.

$\FR{\omega_1}$ has density $\omega$ and $\pi$-weight $\omega_1$, so these may
be different.

We have no information on the possible differences between $t$ and Ind, $\pi_\chi$ and Ind, $\pi$ and cardinality, nor $t$ and $\pi_\chi$.

\chapter{$n$-independence in other Boolean Algebras}
Much as one may study independent sets in \bas{} which are not free, one may
study $n$-independents sets in \bas{} which are not $n$-free.  To this
end we will introduce some cardinal functions.

\section{$\PInd$}
The simplest way to introduce a cardinal function related to
$n$-independence is to measure the size of large $n$-independent
subsets (or, equivalently, $n$-free subalgebras)
 of a \ba{} $A$ with a sup-function.  
\begin{definition}For $1\leq n\leq \omega$, the $n$-independence number of $A$, denoted $\PInd \pa A$, is the supremum of the
cardinalities of the $n$-independent subsets of $A$, that is,
$$\PInd\pa A \defeq \sup\left\{\left|X\right| : X\sub A\; \mbox{and}\; X\;\mbox{is}
\;n\!\mbox{-independent}\right\}$$
\end{definition}

Since every independent set is $n$-independent for $2\leq n\leq \omega$, as is any infinite pairwise
disjoint set (theorem \ref{ipd}), clearly $\PInd\pa A \geq c\left(A\right)+
\mbox{Ind}\left(A\right)$ for infinite $A$. Strict inequality may be possible:  let
$G$ be a graph given by $\omega_1\not\rightarrow\pa{\omega_1}^2_2$; then $G_+$ has only countable independent subsets and only countable pairwise disjoint subsets; the same may hold for $\BC\G$.

 Also, as $m$-independent sets are also $n$-independent for $m\leq n\leq\omega$, $\mbox{mInd}\pa A \leq \mbox{nInd}\pa{A}$.  There is more to say:
\begin{theorem} For any \ba{} $A$, $\PInd\pa A \leq s\pa
A$.\end{theorem}
\begin{proof}

From theorem \ref{pindisidind}, we have that $$\s{\left|X\right|: X\sub A\;\mbox{is}\; n\mbox{-independent}}\sub \s{\left|X\right|: X\sub A\;\mbox{is ideal-independent}},$$ so that their suprema have the relationship $$\sup\s{\left|X\right|: X\sub A\;\mbox{is}\; n\mbox{-independent}}\leq \sup\s{\left|X\right|: X\sub A\;\mbox{is ideal-independent}},$$ that is, $\PInd\pa A \leq s\pa A$.
\end{proof}

We do not know if strict inequality is possible, nor have we investigated the
relationship of $\PInd$ to any other cardinal function.

$\PInd$ behaves reasonably under products.
\begin{theorem}For $A_i$ infinite \bas{}, $\PInd\pa{\prod^{\mathrm w}_{i\in I}A_i}\geq \sum_{i\in I}\PInd\pa{A_i}$\end{theorem}
\begin{proof}
Let $H_i\sub A_i^+$ be $n$-independent.  Then from theorem \ref{bigbindset},
there is an $n$-independent set $\bigcup_{i\in I}H_i$, of size $\sum_{i\in I}\left|H_i\right|$, showing that \\$\PInd\pa{\prod^w_{i\in I}A_i}\geq\sum_{i\in I}\PInd\pa{A_i}$.
\end{proof}

Like Ind, having $\PInd\pa A = \left|A\right|$ does not imply that $A$ is
$n$-free. For $2\leq n\leq \omega$ It is well known that $\P\pa\omega /fin$ has cellularity $2^\omega$, and thus
$\PInd\; 2^\omega$ as well, though it is not $n$-free.

\section{Maximal $n$-independence number}
Since $\PInd$ is a regular sup-function, we can define a spectrum function and a maximal $n$-independence number of a \ba{} in the standard way.
\begin{definition}Let $1\leq n\leq \omega$.
$$\J_{nsp}\pa A \defeq\left\{\left|X\right|: X\;\mbox{is a
maximal}\;n\!\mbox{-independent subset of}\; A\right\}$$
$$\J_n\pa A \defeq \min\pa{\J_{nsp}\pa A}$$
\end{definition}
This could be written as $\PInd_{mm}$ according to the notation of
Monk \cite{CINV2}.  Note that $\J_1=\mathfrak i$ where $\mathfrak i$ is the
minimal independence number as seen in Monk \cite{ContCard}.

This is defined for every \ba{}; from the definition it is easily seen that the
union of a chain of $n$-independent sets is $n$-independent, so Zorn's
lemma shows that there are maximal $n$-independent sets.  It is infinite if $A$ is atomless (shown in lemma \ref{mpi_sum_1}), and has value $1$ if $A$ has an atom.

\begin{proposition} For all $n$ with $1\leq n\leq\omega$, if $A$ has an atom, then $\J_n\pa A=1$.\end{proposition}
\begin{proof}
If $a$ is an atom of $A$, then we claim that $\s{-a}$ is a maximal $n$-independent subset of $A^+$.  That $\s{-a}$ is $n$-independent is clear as any singleton other than $\s 0$ and $\s1$ is independent. 

Let $x\in A^+\sm\s{-a}$, we show that $\s{-a,x}$ is not $n$-independent.
There are two cases.

If $a\leq x$, then $1=a+-a\leq x + -a$, so that $\pa{\perp1}$ fails.

If $a\leq -x$, then $x\leq -a$, so that $0\neq\prod\s x\leq\sum\s{-a}$, but $\s
x\cap\s{-a}=\emptyset$, so that $\pa{\perp3}$ fails.
\end{proof}

\begin{lemma}\label{mpi_sum_1} Let $B$ be a \ba{}, $2\leq n\leq \omega$, and $H\sub B^+$ be
$n$-independent.  If $H$ is maximal among $n$-independent subsets
of $B^+$, then $H$ is infinite and $\sum H =1$ or $H$ is finite and $-\sum H$ is an atom.\end{lemma}
\begin{proof}We prove the contrapositive.
First, the case that $H$ is infinite.
Let $H\sub B^+$ be $n$-independent and have $b<1$ as an upper bound. 
We show that $H\cup\s{-b} $ is $n$-independent:

Note that $-b\notin H$, as $-b\not\leq b$.  Now we will apply proposition
\ref{prop:a}.  So, assume that $R\in\subs{H\cup\s{-b}}{<\omega},
\epsilon\in\leftexp R2$, and $\prod_{x\in R}x^{\epsilon_x}=0$.  If
$-b\notin R$, the conclusion follows since $H$ is $n$-independent.  So suppose
that $-b\in R$.  Let $R'\defeq R\sm\s{-b}$.  Then we have two cases:
\begin{enumerate}
\item[Case 1.] $\epsilon_{-b}=1$.  If there is an $x\in R'$ such that
$\epsilon_x=1$, then $x\leq b$ and so $x\cdot -b =0$ as desired.  So assume
that $\epsilon\left[R'\right]=\s0$.  Then $-b\leq\sum_{x\in R'}x\leq b$, which
is a contradiction.
\item[Case 2.] $\epsilon_{-b}=0$.  If $\epsilon_x=1$ for some $x\in R'$, then 
$$0=\prod_{y\in R} y^{\epsilon_y}=\prod_{y\in R'}y^{\epsilon_y}\cdot
b=\prod_{y\in R'}y^{\epsilon_y}$$
and the $n$-independence of $H$ gives the result.  So assume that
$\epsilon\br{R'}=\s0$.  Then $b\leq\sum R'\leq b$, so $b=\sum R'$.  Then
$z\cdot \prod_{x\in R'}-x=0$, contradicting the $n$-independence of $H$.
\end{enumerate}
So we have that if $H$ is infinite and maximal $n$-independent, it has no upper bound
other than $1$, so $\sum H =1$.

Now we consider the case that $H$ is finite.
 If $-\sum H$ is not an atom, let $0<a<-\sum H$, then we claim that $H\cup\s a$
is $n$-independent. Again we use proposition \ref{prop:a}.  Assume that
$R\in\subs{H\cup\s a}{<\omega}, \epsilon\in\leftexp R2$, and $\prod_{x\in
R}x^{\epsilon_x}=0$. \Wolog, $a\in R$. 
\begin{enumerate}
\item[Case 1.]$\epsilon_a=1$.  If $\epsilon_x=1$ for some $x\in R\sm\s a$, then
$a\cdot x\leq a\cdot\sum H=0$, as desired.  Otherwise $a\leq \sum\pa{R\sm\s
a}\leq \sum H$ and so $a=0$, contradiction.
\item[Case 2.]$\epsilon_a=0$.  If $\epsilon_x=1$ for some $x\in R\sm\s a$, then
$a\cdot x=0$, hence $x\leq-a$, and then $\prod_{y\in
R}y^{\epsilon_y}=\prod\s{y^{\epsilon_y}:y\in R\sm\s a}$ and the conclusion
follows.  Otherwise $-a\leq\sum\pa{R\sm\s a}\leq\sum H$, so $-\sum H\leq a$,
contradicting $a<-\sum H$.
\end{enumerate}

\end{proof}

The converse of lemma \ref{mpi_sum_1} does not hold.  An example due to Monk is in  $\FR{X\cup Y}$ where $X\cap Y=\emptyset$ and 
 $\left|X\right|=\left|Y\right|=\kappa\geq\omega$.
  $X$ is independent, is not maximal for $2$-independence, and has sum 1.  Here $\sum X=1$ is the only non-trivial part--by way of contradiction, let $b$ be a non-1 upper bound for $X$.  Then $-b$ has the property that $x\cdot -b =0$ for all $x\in X$, so let $a$ be a elementary product of elements of $X\cup Y$ where $a\leq -b$.  Take some $x\in X$ that does not occur in that elementary product.  Then since $X\cup Y$ is independent, $a\cdot x \neq 0$, but since $a\leq -b$, $a\cdot x =0$.

\begin{theorem}For $B$ atomless, and $2\leq n\leq \omega$, $\mathfrak p\pa B \leq \J_n \pa B$.\end{theorem}
Here $\mathfrak p\pa B$ is the pseudo-intersection number, defined in
Monk \cite{ContCard} as 
$$\mathfrak p \pa A \defeq \min\s{\left|Y\right|:\sum Y = 1 \;\mbox{and}\; \sum Y'\neq 1\;
\mbox{for every finite}\; Y'\sub Y}.$$
\begin{proof}
Since $B$ is atomless, a maximal $n$-independent set $Y$ has $\sum Y=1$,
and by \pc1, if $Y'\sub Y$ is finite, $\sum Y'\neq 1$.  That is, the
maximal $n$-independent sets are included among the $Y$ in the definition
of $\mathfrak p\pa A$. 
\end{proof}

We do not know if strict inequality is possible.

\begin{corollary} For all $n$ with $1\leq n\leq \omega$, $\J_n\pa{\P\pa\omega/fin}\geq\aleph_1$\end{corollary}
\begin{proof}$\aleph_1\leq \mathfrak p \pa{\P\pa\omega/fin} \leq \J_n\pa{\P\pa\omega/fin}$
\end{proof}

We also recall that under Martin's Axiom, $\mathfrak p \pa{\P\pa\omega/fin} =
\beth_1$, so the same is true of $\J_n$.

\begin{proposition} Any $B$ with 
the strong \CSP{} has, for all $2\leq n\leq \omega$,  $\J_n\pa B\geq \aleph_1$.\end{proposition}
\begin{proof}Such a $B$ is atomless, so let $H\sub B^+$ be
$n$-independent and countably infinite, that is $H=\left<h_i:i\in \omega\right>$. Then
let $c_m\defeq\sum_{i\leq m} h_i$. Each $c_m$ is a finite sum of elements of $H$,
thus by $\pa{\perp1}$, $c_m<1$. Then $C\defeq\s{c_i:i\in \omega}$ is a countable
chain in $B\sm \s1$, so by the strong \CSP{}, there is a $b\in B$ such that
$c_i\leq b <1$ for all $i\in \omega$.  Then as $h_i\leq c_i$, $h_i\leq b$ for
all $i\in\omega$ as well, that is, $b$ is an upper bound for $H$.  Thus by
lemma \ref{mpi_sum_1}, $H$ is not maximal.
\end{proof}

In addition, we show that maximal $n$-independent sets lead to weakly dense sets. 

We use the notation $-X=\s{-x : x\in X}$ frequently in the sequel.

\begin{theorem}Let $1\leq n\leq\omega$.  If $X\sub A$ is maximal $n$-independent in $A$, then the set $Y$ of nonzero elementary products of elements of $X$ is weakly dense in $A$.
\end{theorem}
Recall that $Y$ is weakly dense in $A$ if and only if $Y\sub A^+$ and for every $a\in A^+$, there is a $y\in Y$ such that $y\leq a$ or $y\leq -a$.
\begin{proof}
If $a\in X$, this is trivial, so we may assume that $a\notin X$ and hence $X\cup\s a$ is not $n$-independent.

By proposition \ref{prop:a}, there exist $R\in\subs{X\cup\s a}{<\omega}$ and
$\epsilon\in \leftexp R2$ such that $\prod_{x\in R}x^{\epsilon_x}=0$ while for
every $R'\in\subs R{\leq n}$, if $\epsilon\br{R'}\sub \s 1$ then $\prod
R'\neq0$.  This last implication holds for every $R'\in\subs{R\sm\s a}{\leq
n}$, and so $\prod\s{x^{\epsilon_x}:x\in R\sm\s a}\neq 0$ since $X$ is
$n$-independent.  But $\prod\s{x^{\epsilon_x}:x\in R\sm\s a}\leq a$ or $\leq
-a$, as desired.
\end{proof}

\begin{corollary}
If $A$ is atomless and $1\leq n\leq \omega$, then $\mathfrak r\pa A\leq \J_n\pa A$.
\end{corollary}
Recall the definition of the reaping number: $$\mathfrak r\pa A\defeq\min\s{\left|X\right| : X\mbox{ is weakly dense in } A}.$$
\begin{proof}
Since $A$ is atomless, all maximal $n$-independent sets are infinite, and thus there is a set of size $\J_n\pa A$ weakly dense in $A$.
\end{proof}
We do not know if strict inequality is possible.

We do not currently have any results for the behavior of $\J_n$ on any type of
product or its relationship to $\mathfrak u$.

We show the consistency of $\J_k\pa{\P\pa \omega /fin} < \beth_1$ for $1\leq k\leq\omega$.  The argument is similar to
exercises (A12) and (A13) in chapter VIII of Kunen \cite{KunenSet}; the main lemma
follows.
 
\begin{lemma}\label{con_j_lemma}
Let $M$ be a countable transitive model of $ZFC$ and $1\leq k\leq\omega$.  For a subset $a$ of
$\omega$, let $\left[a\right]$ denote its equivalence class in $\P\pa\omega /
fin$.  Suppose that $\kappa$ is an
infinite cardinal and $\left<a_i : i<\kappa\right>$ is a system of infinite
subsets of $\omega$ such that $\left<\left[a_i\right]:i<\kappa\right>$ is
$k$-independent in $\P\pa\omega / fin$.  Then there is a generic extension
$M\left[G\right]$ of $M$ using a ccc partial order such that in
$M\left[G\right]$ there is a $d\sub \omega$ with the following properties:
\begin{enumerate}
\item $\left<\left[a_i\right]:i<\kappa\right>^\frown \left<\left[\omega\sm
d\right]\right>$ is $k$-independent.

\item If $x\in \pa{\P\pa\omega \cap M}\sm\pa{\s{a_i : i<\kappa}\cup\s{\omega\sm
d}}$, then $\left<\left[a_i\right]:i<\kappa \right>^\frown\left<\left[\omega\sm
d\right],\left[x\right]\right>$ is not $k$-independent.
\end{enumerate} 
\end{lemma}
\begin{proof}We work within $M$ here.
Let $B$ be the $k$-independent subalgebra of\\ $\P\pa\omega / fin$ generated by $\s{\left[a_i\right]:i<\kappa}$.  By Sikorski's extension criterion, let $f$ be a homomorphism from $\left<\s{a_i : i<\kappa}\cup \s{\s m : m\in \omega}\right>$ to $\overline B$ such that $f\pa{a_i}=\left[a_i\right]$ and $f\pa{\s m}=0$.  
Then let $h:\P\pa\omega \fct \overline B$ be a homomorphic extension of $f$ as given by Sikorski's extension theorem.

Let $ P\defeq\s{\pa{b,y}:b\in\ker\pa h \mbox{ and }
y\in\left[\omega\right]^{<\omega}}$ with the partial order given by
$\pa{b,y}\leq\pa{b',y'}$ if and only if $b\supseteq b'$, $y\supseteq y'$ and $y\cap b'\sub
y'$.  This is a ccc partial order.  Let $G$ be a $ P$-generic filter
over $M$, and let $d\defeq\bigcup_{\pa{b,y}\in G}y$.

We now have several claims that combine to prove the lemma.

\begin{itemize}

\item[Claim 1.] If $R$ is a finite subset of $\kappa$ and $\epsilon\in \leftexp R 2$ is such that $\bigcap_{\stackrel{i\in R}{\epsilon_i=1}}a_i$ is infinite, then $\bigcap_{i\in R}a_i^{\epsilon_i}\cap d$ is infinite.

Let $R$ and $\epsilon$ be as given, then for each $n\in \omega$, let 
$$E_n\defeq \s{\pa{b,y}\in P : \exists m>n\left[m\in\bigcap_{i\in R}a_i^{\epsilon_i}\cap y\right]}.$$

First, we show that each $E_n$ is dense.  Take $\pa{b,y}\in  P$.  Then
$c\defeq\pa{\bigcap_{i\in R}\pa{a_i^{\epsilon_i}}}\sm b$ is infinite; if
not, then $c$ is finite (thus in $\ker\pa h$, as is $b$)and $\bigcap_{i\in R}a_i^{\epsilon_i}\sub b\cup c$.  Applying $h$ to both sides gives $\prod_{i\in R}\left[a_i\right]^{\epsilon_i}=0$, which is a contradiction of proposition \ref{prop:a}.  So we choose an $m\in c\sm y$ such that $m>n$; then $\pa{b,y\cup\s m}\leq\pa{b,y}$ and  $\pa{b,y\cup\s m}\in E_n$, showing that $E_n$ is dense. 
This shows the claim, as for each $n\in\omega$, $E_n\cap G\neq \emptyset$, so that we have an integer larger than $n$ in $\bigcap_{i\in R}a_i\cap d$.

\item[Claim 2.]  If $R$ is a finite subset of $\kappa$ and $\epsilon\in \leftexp R 2$ such that $\bigcap_{\stackrel{i\in R}{\epsilon_i=1}}a_i$ is infinite, then $\bigcap_{i\in R}a_i^{\epsilon_i}\sm d$ is infinite. 

Let $R$ and $\epsilon$ be as given, then for each $n\in \omega$, let 
$$D_n \defeq\s{\pa{b,y}\in P:\exists m>n\left[m\in\bigcap_{i\in R}a_i^{\epsilon_i}\cap b\sm y\right]}. $$

To show that $D_n$ is dense, take any $\pa{b,y}\in P$.  Since $\bigcap_{i\in
R}a_i^{\epsilon_i}$ is infinite from proposition \ref{prop:a}, it follows that
we may choose $m>n$ such that $m\in \bigcap_{i\in R}a_i^{\epsilon_i}\sm y$.
Then $\pa{b\cup\s m , y}\leq\pa{b,y}$ and $\pa{b\cup\s m,y}\in D_n$, as
desired.

Take some $\pa{b,y}\in D_n\cap G$.  Then there is an $m>n$ such that $m\notin d$ (thus proving the claim). In fact, choose $m>n$ such that $m\in\bigcap_{i\in R}a_i^{\epsilon_i}\cap b \sm y$. We claim that $m\not\in d$. Suppose that $m\in d$;
 then we have a $\pa{c,z}\in G$ with $m\in z$ and $\pa{e,w}\in G$ that is a common extension of $\pa{b,y}$ and $\pa{c,z}$.  Then $m\in w\cap b\sm y$, contradicting that $\pa{e,w}\leq\pa{b,y}$.

\item[Claim 3.]
$\left<\left[a_i\right]:i<\kappa\right>^\frown\left<\left[\omega\sm
d\right]\right>$ is $k$-independent.

Suppose that $R\in\subs\kappa{<\omega},\epsilon\in\leftexp R2, \delta\in 2$,
and $\prod_{i\in R}\br{a_i}^{\epsilon_i}\cdot\br{\omega\sm d}^\delta =0$.  By
claims 1 and 2 (depending on $\delta$), $\prod_{\stackrel{i\in
R}{\epsilon_i=1}}\br{a_i}^{\epsilon_i}=0$.  Since
$\left<\br{a_i}:i<\kappa\right>$ is $k$-independent, there is a subset
$R'\sub\s{i\in R:\epsilon_i=1}$ of size at most $k$ such that $\prod_{i\in
R'}\br{a_i}=0$, as desired.

\item[Claim 4.] If $b\in \ker\pa h$, then $b\cap d$ is finite.

$\s{\pa{c,y}\in P:b\sub c}$ is dense in $P$, so that there is a $\pa{c,y}\in G$
such that $b\sub c$.  We show $b\cap d \sub y$ and thus is finite.  Let $m\in
b\cap d$ and choose an $\pa{e,z}\in G$ such that $m\in z$.  Let $\pa{r,w}\in G$
be a common extension of $\pa{e,z}$ and $\pa{c,y}$; then (recalling the
definition of the order) $m\in w\cap c\sub y$.

\item[Claim 5.]If $x\in
\pa{\P\pa\omega\cap M}\sm\pa{\s{a_i:i<\kappa}\cup\s{\omega\sm d}}$, then\\
$s\defeq\left<\left[a_i\right]:i<\kappa\right>^\frown\left<\left[\omega\sm
d\right],\left[x\right]\right>$ is not $k$-independent.

We have two cases here.  The slightly easier is if $x\in\ker\pa h$; then by claim 4,
$x\cap d$ is finite, so that $\left[x\right]\leq\left[\omega\sm d\right]$,
 causing $s$ to fail to even be ideal-independent.  If
$x\notin\ker\pa h$, then there is a $b\in B$ with
$0<b\leq h\pa x$.  Since $B$ is $k$-freely generated by
$\left<\left[a_i\right]:i<\kappa\right>$, we may take $b$ to be a elementary product
of elements of $\left<\left[a_i\right]:i<\kappa\right>$.  Then
$b=\left[c\right]$, where $c=\bigcap_{i\in R}a_i^{\epsilon_i}$ is infinite.  Then $c\sm x\in \ker\pa h$.  By claim 4, this gives $\prod_{i\in R}\left[a_i\right]^{\epsilon_i}\cdot -\left[x\right]\cdot\left[d\right]=0$, contradicting proposition \ref{prop:a} for $s$.

\end{itemize}
\end{proof}

\begin{theorem}For each $1\leq k\leq\omega$, it is consistent with $\beth_1>\aleph_1$ that $\J_k\pa{\P\pa\omega/fin}=\aleph_1$.\end{theorem}
\begin{proof}We begin with a countable transitive model $M$ of $ZFC + \beth_1>\aleph_1$, then iterate the construction of lemma \ref{con_j_lemma} $\omega_1$ times as in lemma 5.14 of chapter VIII of Kunen \cite{KunenSet}. This results in a model of $ZFC + \beth_1>\aleph_1 + \J_k\pa{\P\pa\omega/fin}=\aleph_1$.
\end{proof}

This shows that $\J_k\pa{\P\pa\omega/fin}=\beth_1$ is independent of $ZFC$.

\chapter{Ideal Independence}
\section{A necessary condition for maximality}
We give a solution to a problem posed in a draft of Monk \cite{minimax}.

\begin{theorem}\label{thm_idindsum1}Let $A$ be a \ba{}. If $X\sub A$ is maximal for
ideal-independence, then $\sum X=1$.
\end{theorem}
\begin{proof}
Let $X\sub A$ be ideal-independent and let $b\in A\sm\s 1$ be such that
$\forall x\in X x\leq b$.  We claim that $X'\defeq X\cup\s{-b}$ is ideal-independent
(thus proving the theorem).

We need to show that for distinct $x,x_1,x_2,\ldots,x_n\in X'$, $x\not\leq
x_1+x_2+\ldots+x_n$.  There are three cases to consider.

\begin{description}
\item{Case I} $-b\not\in\s{x,x_1,x_2,\ldots,x_n}$.

Then $x\not\leq x_1+x_2+\ldots+x_n$ by the ideal-independence of $X$.

\item{Case II} $-b=x$.

Assume otherwise, that is, $-b\leq x_1+x_2+\ldots+x_n$.
Since $x_1+x_2+\ldots+x_n\leq b$, we have $-b\leq b$, thus $b=1$,
contradiction.

\item{Case III} $-b=x_1$.

Then as $x\in X$, $x\leq b$, that is, $x\cdot -b = 0$.
By the ideal-independence of $X$, $x\not\leq x_2 + x_3 + \ldots + x_n$, so that
$y\defeq x\cdot\pa{x_2 + x_3 + \ldots + x_n}<x$.

Then $x\cdot\pa{-b + x_2 + x_3 + \ldots + x_n}=x\cdot -b + x\cdot \pa{x_2 + x_3
+ \ldots + x_n}=0+y<x$, thus $x\not\leq -b + x_2 + x_3 + \ldots + x_n$.
\end{description}
\end{proof}

Recall that $$s_{mm}\pa A \defeq\min\s{\left|X\right|: X\mbox{ is
maximal ideal-independent}}.$$

\begin{corollary}
If $A$ is infinite, then $s_{mm}\pa A$ is infinite.
\end{corollary}

\begin{corollary}
For $A$ infinite, $\mathfrak p\pa A \leq s_{mm}\pa A$
\end{corollary}
Recall that
$$\mathfrak p \pa A \defeq \min\s{\left|Y\right|:\sum Y = 1 \;\mbox{and}\; \sum Y'\neq 1\;
\mbox{for every finite}\; Y'\sub Y}.$$

\begin{proof}
Let $X$ be maximal ideal-independent (and thus infinite).  By theorem \ref{thm_idindsum1}, $\sum
X=1$.  Also, if $X'\sub X$ is finite, $\sum X'\neq 1$, otherwise, take $x\in
X\sm X'$, then $x\leq 1 = \sum X'$, contradicting the ideal-independence of
$X$.  Thus $$\s{\left|X\right|: X\mbox{ is
maximal ideal-independent}}\sub$$ $$\s{\left|Y\right|:\sum Y = 1 \;\mbox{and}\; \sum Y'\neq 1\;
\mbox{for every finite}\; Y'\sub Y}$$ and so $\mathfrak p\pa A \leq s_{mm}\pa
A.$
\end{proof}
We do not know of a relation between $\J_n$ and $s_{mm}$.

\chapter{Moderation in Boolean Algebras}
\section{Definitions}
Moderateness was explicitly defined by Heindorf \cite{Moderate} to systematize an often used construction.

\begin{definition} Let $A$ be a \ba{}, and $F,G\sub A$.
\begin{enumerate}
\item For any $a\in A$, its norm with respect to $F$, $\nor a_F$ is defined as $$\nor a_F \defeq\s{f\in F : 0<a\cdot f < f}.$$  The subscript may be dropped if $F$ is clear in context.  When $f\in \nor a_F$, it is said that ``$a$ splits $f$'' since $f\cdot -a\neq 0$ and $f\cdot a\neq 0$.

\item $F$ is moderate in $G$ if and only if $\nor g_F$ is finite for every $g\in G$.

\item $F$ is moderate if and only if $F$ is moderate in $A$ and $0\notin F$.

\item $a\in A$ is saturated with respect to $F$ if and only if $\nor a_F=\emptyset$.
\end{enumerate}
\end{definition}

The following obvious but useful lemma and corollary are 1.3 and 1.4 in Heindorf \cite{Moderate}.

\begin{lemma}$$\nor a = \nor{-a},\;\; \nor{a+b}\sub \nor a \cup \nor b, \mbox{and} \;\; \nor{a\cdot b}\sub \nor a \cup \nor b$$
\end{lemma}

\begin{corollary}\label{cor:initself}$F$ is moderate in $G$ if and only if $F$ is moderate in $\left<G\right>$.\end{corollary}

We also use the shorthand ``$F$ is moderate in itself'' for $F$ is moderate in
$F$, as well as the elementary fact that if $F$ is moderate in itself and
$G\sub F$, then $G$ is moderate in itself.

The construction inspiring these definitions is named the moderate product, and is related to moderate algebras, defined as those algebras in which every ideal is generated by a moderate set.

We will consider ``moderately generated algebras'', that is, those algebras which are generated by a moderate set.  By corollary \ref{cor:initself}, these are the algebras generated by sets which are moderate in themselves.  Examples of such sets are relatively easy to come by.  For example, any disjoint set is moderate in itself, so finite-cofinite algebras are moderately generated.  Countable \bas{} are also included in this class:

\begin{proposition}Let $A$ be a countable \ba{}; then $A$ is moderately generated.
\end{proposition}
\begin{proof}Let $X$ be a generating set for $A$.  Fix an enumeration of $X=\left<x_i:i<\omega\right>$, and let $$F\defeq\s{\prod_{i\leq n}x_i^{\epsilon_i} : n\in\omega, \epsilon\in\leftexp n2}.$$

We claim that $F$ generates $A$.  We prove this claim by showing that $X\sub \left<F\right>$.

Clearly $x_0\in F\sub\left<F\right>$.  For $n\geq 1$, $x_n=\sum_{\epsilon\in\leftexp n 2}\pa{\prod_{i<n}\pa{x_i^{\epsilon_i} \cdot x_n}}.$  The products in the right hand side are all in $F$, hence the total is in $\left<F\right>$.

We claim that $F$ is moderate in $F$.  We prove this claim by calculating $\nor{x_m}_F$.

$$\nor{x_0}_F = \s{\prod_i x_i^{\epsilon_i} : 0 < x_0 \cdot \prod_i x_i^{\epsilon_i} < \prod_i x_i^{\epsilon_i}}$$  This set is empty--if $\epsilon_0=1$, then $ x_0 \cdot \prod_i x_i^{\epsilon_i} = \prod_i x_i^{\epsilon_i}$ and if $\epsilon_0=0$, then $0= x_0 \cdot \prod_i x_i^{\epsilon_i}$.

For $m\geq 1$, we claim that 
$$\nor{x_m}_F\sub\s{\prod_{i\leq n} x_i^{\epsilon_i} : n<m, \epsilon\in \leftexp n2},$$
which is of size $\leq 2^{m+1}$. We take $n\geq m$ and $\epsilon\in \leftexp n2$ and show that $\prod_{i\leq n}x_i^{\epsilon_i}\notin \nor{x_m}_F$.  Since $n\geq m$, then we may consider values of $\epsilon_m$. If $\epsilon_m=1$, then $x_m\cdot\prod_{i\leq n} x_i^{\epsilon_i} = \prod_{i\leq n} x_i^{\epsilon_i}$.  If $\epsilon_m=0$, then $x_m\cdot\prod_{i\leq n} x_i^{\epsilon_i}=x_m\cdot x_m^0 \cdot \prod_{i\leq n} x_i^{\epsilon_i}=0$.  Thus $\prod_{i\leq n}x_i^{\epsilon_i}\notin \nor{x_m}_F$.
\end{proof}

An infinite independent set is not moderate in itself; in fact if $X$ is infinite independent and
$x\in X$, then $\nor x_X=X\sm\s x$.  Note that $x\notin \nor x$, so this is as
large as $\nor x$ can ever be.  Thus $\FR\kappa$ for $\kappa\geq\aleph_1$
is not moderately generated since any generating set must be of size $\kappa$
and thus by theorem 9.16 of Koppelberg \cite{Handbook}, must include an
independent set of size $\kappa$.

A $2$-independent set has a better chance of being moderate in itself (of
course a disjoint set is $2$-independent and each element of such a set is
saturated with respect to the disjoint set), as it is again simple to compute
norms.  Let $X$ be $2$-independent; then for $x\in X$, since $X$ is
incomparable, $\nor x_X=\s{y\in X : x\cdot y\neq 0}$.  That is, if
$\left<X\right>=\BC\G$, and $v$ is a vertex of $\G$ with neighborhood $N\pa v$, then $\nor{v_+}_X =
N\pa{v}_+$.  So a $2$-independent set is moderate exactly when its intersection graph is such that every vertex has finite degree.

\section{Incomparability}

In the previous section, we noted that a pairwise disjoint set is moderate in
itself.   We will prove a partial converse: most moderate sets have a large
incomparable subset.

\begin{proposition}Let $F$ be  moderate in itself in a \ba{} $A$ such that $\left|F\right|$ is an uncountable regular cardinal.  Then $F$ has an incomparable subset of size $\left|F\right|$.
\end{proposition}
\begin{proof} Let $f\in F$.  Since $\nor{f}_F$ is finite, the set $X_f\defeq\s{a\in F:f\leq a}\sub\s f\cup\nor{f}_F$ is finite as well.  We first show that if $f\neq g$, then $X_f\neq X_g$:  If $X_f\sub X_g$, then for all $a\in F$, $f\leq a$ implies that $g\leq a$; in particular, $g\leq f$.  Thus if $X_f=X_g$, then $f=g$.
So the collection $\mathcal X\defeq\s{X_f:f\in F}$ has the same cardinality as $F$.  Thus, it has a subset $\mathcal Y$ which is a $\Delta$-system with root $r$.  Let $H\defeq \bigcup\s{X_f\sm r : X_f\in \mathcal Y}$.  
$\left|H\right|=\left|F\right|=\left|\mathcal Y\right|$.  
Let $Y_f\defeq X_f\sm r$ for each $X_f\in \mathcal Y$ and note that $Y_f\cap Y_g=\emptyset$ for distinct $f,g$. Since $r$ is finite, without loss of generality, we may assume that $f\in Y_f$. 
In particular, $f\notin Y_g$ and $g\notin Y_f$, that is, $g\not\leq f$ and $f\not\leq g$, so that $G\defeq\s{f:X_f\in \mathcal Y}$ is an incomparable set of size $\left| F\right|.$

\end{proof}

\begin{corollary}If $\kappa$ is the cardinality of an uncountable moderate family of regular cardinality in $A$, then $\kappa\leq \mbox{Inc}\pa A$.\end{corollary}

\begin{corollary}If $A$ is moderately generated, then $\mbox{Inc}\pa A =
\left|A\right|$.  Incomparability is attained if $\left|A\right|$ is regular.\end{corollary}
From this, we may also conclude that
h-cof is equal to cardinality for moderately generated \bas{}.

\bibliographystyle{plain}
\bibliography{biblio}
\end{document}